%% file: main.tex
\documentclass[10pt]{article} 
\usepackage[preprint]{tmlr}

\title{Solving Imaging Inverse Problems Using Plug-and-Play Denoisers: Regularization and Optimization Perspectives}

\author{\name Hong Ye Tan \email hyt35@cam.ac.uk \\
      \addr 
      University of Cambridge
      \AND
      \name Subhadip Mukherjee \email smukherjee@ece.iitkgp.ac.in \\
      \addr IIT Kharagpur
      \AND
      \name Junqi Tang \email j.tang.2@bham.ac.uk \\
      \addr University of Birmingham}

\usepackage{microtype}
\usepackage{amsmath,amssymb}
\usepackage{mathtools}
\usepackage{mathrsfs}
\usepackage{xcolor}
\usepackage{tikz}
\usepackage{mathrsfs}

\usepackage{graphicx}
\usepackage{natbib}
\usepackage{enumitem}
\usepackage{algorithm}
\usepackage{algpseudocode}
\usepackage{hyperref}
\usepackage{cleveref}
\usepackage{subcaption}
\usepackage[normalem]{ulem}
\usepackage{url}
\usepackage{pifont}
\usepackage{multirow}
\usepackage{booktabs}
\newcommand{\cmark}{\ding{51}}%
\newcommand{\xmark}{\ding{55}}%
\newcommand{\R}{\mathbb{R}}
\DeclareMathOperator*{\argmin}{arg\,min}

\DeclareMathOperator*{\zer}{zer}
\DeclareMathOperator*{\dom}{dom}
\DeclareMathOperator{\prox}{prox}
\DeclareMathOperator{\dist}{dist}
\DeclareMathOperator{\id}{Id}
\DeclareMathOperator{\Fix}{Fix}
\DeclareMathOperator{\PnP}{PnP}

\newcommand{\Op}[1]{\operatorname{\mathcal{#1}}}

\newtheorem{theorem}{Theorem}[section]
\newtheorem{definition}[theorem]{Definition}
\newtheorem{proposition}[theorem]{Proposition}
\newtheorem{lemma}[theorem]{Lemma}

\newtheorem{remark}[theorem]{Remark}
\newcommand{\rev}[1]{#1}
\begin{document}

\maketitle

\begin{abstract}
Inverse problems lie at the heart of modern imaging science, with broad applications in areas such as medical imaging, remote sensing, and microscopy. Recent years have witnessed a paradigm shift in solving imaging inverse problems, where data-driven regularizers are used increasingly, leading to remarkably high-fidelity reconstruction. A particularly notable approach for data-driven regularization is to use learned image denoisers as implicit priors in iterative image reconstruction algorithms. This chapter presents a comprehensive overview of this powerful and emerging class of algorithms, commonly referred to as plug-and-play (PnP) methods. We begin by providing a brief background on image denoising and inverse problems, followed by a short review of traditional regularization strategies. We then explore how proximal splitting algorithms, such as the alternating direction method of multipliers (ADMM) and proximal gradient descent (PGD), can naturally accommodate learned denoisers in place of proximal operators, and under what conditions such replacements preserve convergence. The role of Tweedie's formula in connecting optimal Gaussian denoisers and score estimation is discussed, which lays the foundation for regularization-by-denoising (RED) and more recent diffusion-based posterior sampling methods. We discuss theoretical advances regarding the convergence of PnP algorithms, both within the RED and proximal settings, emphasizing the structural assumptions that the denoiser must satisfy for convergence, such as non-expansiveness, Lipschitz continuity, and local homogeneity. We also address practical considerations in algorithm design, including choices of denoiser architecture and acceleration strategies. By integrating both classical optimization insights and modern learning-based priors, this chapter aims to provide a unified and accessible framework for understanding PnP methods and their theoretical underpinnings.  
\end{abstract}

\include{chapter1_new}

\end{document}

%% file: chapter1_new.tex
\section{Introduction}
Inverse problems are the backbone of numerous applications in imaging science, signal processing, computational physics, and beyond. In a typical \textit{ill-posed} imaging inverse problem, one seeks to recover an unknown signal or image $x$ from its indirect and often noisy observations $y = Kx+w$, where $K$ denotes the forward operator representing the imaging process and $w$ denotes measurement noise. The image $x$ and its measurement $y$ are assumed to lie in appropriate normed vector spaces. Classic examples include low-level computer vision tasks, such as image deblurring, denoising, super-resolution, and inpainting, as well as various medical imaging tasks, including MRI, PET, and CT reconstruction. Due to incomplete or corrupted measurements, inverse problems are fundamentally ill-posed and require \textit{regularization} to ensure stable and meaningful solutions.

Traditional approaches to regularization (see \cite{engl1996regularization,scherzer2009variational,benning2018modern} for a detailed treatment) rely on the design of explicit prior models that encode different regularity assumptions on images, such as Tikhonov regularization, total variation (TV), or, more recently, sparsity-promoting regularizers (typically penalizing the $\ell_1$-norm in a suitable transform domain). These classical regularization approaches, while theoretically grounded and computationally tractable, often fail to fully capture the rich and highly structured nature of natural images, thereby limiting reconstruction quality in challenging scenarios.

Thanks to the pioneering work by Venkatakrishnan et al. \cite{venkat_pnp_6737048}, \emph{plug-and-play} (PnP) methods have emerged in the last decade as an effective paradigm, deviating from handcrafted priors by instead integrating powerful image denoisers directly into iterative optimization algorithms for solving inverse problems (see \cite{kamilov2023plug} for a recent survey on PnP methods). The key insight behind PnP is deceptively simple yet quite powerful: instead of specifying an explicit prior, one can ``plug in'' an off-the-shelf image denoiser within an iterative framework such as \textit{alternating direction method of multipliers} (ADMM) or \textit{proximal gradient descent} (PGD), effectively regularizing the solution through repeated denoising operations. This idea has led to a family of flexible algorithms, with PnP-ADMM~\cite{chan2016plug} and PnP-PGD \cite{pmlr-v97-ryu19a,convergent_prox_pnp_pgd} being two of the most prominent and widely studied variants of PnP algorithms. PnP methods are particularly attractive due to their modularity; the same image denoiser leads to a reasonable image reconstruction for different forward operators. By decoupling the forward model from the prior, PnP enables practitioners to exploit state-of-the-art denoisers, whether model-based (such as BM3D \cite{bm3d_main}) or data-driven (e.g. deep convolutional neural network (CNN)-based denoisers \cite{dncnn_main,pnp_deep_denoiser}), without the need to derive complex closed-form priors or incorporate explicit regularizers into the variational framework. This practical plug-and-play nature makes these methods highly attractive and flexible for diverse applications where noise statistics, measurement models, and image features vary widely.

In parallel, the \emph{regularization by denoising} (RED) approach ~\cite{romano2017little,red_schniter} offers another compelling framework by constructing a class of explicit regularizers using off-the-shelf image denoisers. RED provides a unifying view that connects variational principles with denoising-based operators, offering fixed-point iterations with interpretability and, under suitable conditions, convergence guarantees~\cite{red_schniter}. Other notable theoretical advances include rigorous convergence guarantees for PnP-PGD and PnP-ADMM under appropriate Lipschitz continuity assumptions on the denoising residual (see Theorems 1, 2, and Corollary 3 in \cite{pmlr-v97-ryu19a}), and the extension of PnP concepts to broader classes of iterative schemes, such as consensus equilibrium \cite{buzzard2018plug} and block coordinate methods~\cite{red_block_coord}. The synergy between PnP frameworks and deep learning has further amplified their impact. Learned denoisers, such as DnCNN~\cite{dncnn_main}, FFDNet \cite{ffdnet}, and more recent transformer-based \cite{pnp_transformer} or denoising diffusion model-based approaches \cite{pnp_denoising_diffusion} can be seamlessly integrated into PnP schemes, delivering remarkable reconstruction performance even in severely ill-posed settings. Beyond static images, PnP methods are increasingly being adapted to dynamic and multimodal imaging tasks \cite{pnp_dynamic_NLOS}, where additional temporal or cross-modality constraints may be incorporated, while still retaining the modular plug-and-play philosophy.

More recently, the emergence of generative models, particularly diffusion probabilistic models, has opened new directions for solving inverse problems (see \cite{daras2024surveydiffusionmodelsinverse} for a recent survey, and \cite{Chung2025} for applications in medical imaging), especially using PnP-like approaches. These models provide powerful priors that can sample high-quality images conditioned on measurements, blending generative sampling with the PnP idea to tackle complex, high-dimensional inverse problems with improved uncertainty quantification.

Despite their empirical success, PnP methods still present open challenges. Key questions include understanding their theoretical guarantees when using highly nonlinear or non-expansive denoisers, designing adaptive schemes that can select or learn denoising strength on the fly, and extending the framework to handle non-Gaussian noise, physics-based constraints, or multimodal data fusion. In this chapter, we aim to provide a comprehensive and up-to-date account of plug-and-play methods for imaging inverse problems. We trace the historical development of this framework, highlight the underlying theoretical foundations, and discuss advances in algorithmic design and learning-based denoisers. We also discuss the application of PnP denoisers for Bayesian imaging that incorporates denoisers with diffusion-based posterior sampling, and identify open research directions that may shape the future of PnP-based inverse problem solving. We focus particularly on methods that are widely regarded as pioneering in the area of PnP imaging and methods that come with rigorous convergence guarantees.
\subsection{Brief Survey of Image Denoising}
\label{sec_brief_denoising_survey}
Image denoising has long been a fundamental problem in signal and image processing, driving the development of increasingly sophisticated models and algorithms over the past several decades. We refer interested readers to \cite{GOYAL2020220_denoise_elsv,elad_denoising_review,is_denoising_dead} for a review on image denoisers, covering classical techniques to modern deep learning-based approaches. In the discussion that follows, we will assume that the image $x$ is discretized, and can therefore be represented as a vector in $\mathbb{R}^n$ \rev{after concatenating all the pixels and the color channels in the image}. \rev{For a grayscale image of size $n_1\times n_2$, $n=n_1 n_2$; and for a color image (with three color channels) of the same size, $n=3n_1n_2$.} Under an additive noise model, the goal of image denoising is to recover an unknown clean image \( x \in \mathbb{R}^n \) from a noisy observation \( z = x + w \), where \( w \) represents additive noise. Depending on the application, $w$ may follow a Gaussian or non-Gaussian distribution, and may be white or colored. The fundamental challenge arises from the fact that noise corrupts both low- and high-frequency components in the image, making it difficult to separate noise from fine image details. Designing an effective denoising algorithm therefore requires balancing noise suppression with the preservation of edges, textures, and small structures in the underlying image. An important feature of well-designed (though not necessarily perfect) denoisers is that they can naturally generate a multiscale decomposition of an image \cite{Milanfar2025-xp}, while still allowing exact reconstruction. This idea parallels classical multiscale decompositions, such as the Laplacian pyramid \cite{laplace_pyramid}, but is achieved here through the action of the denoiser.

The earliest and perhaps most intuitive family of denoising algorithms is based on linear filtering (in spatial or frequency domains) \cite[Chapters 3, 4]{gonzalez2008digital}. Classical linear filters, such as the Gaussian filter, exploit the assumption that noise essentially has largely high-frequency components, while natural images exhibit local spatial smoothness (which represents low-frequency features). These filters convolve the noisy image with a spatially localized kernel, attenuating high-frequency components. Mathematically, the output can be expressed as \( \hat{x} = H z \), where \( H \) denotes the convolution operator defined by a filter kernel. Although simple and computationally efficient, linear smoothing tends to blur edges and oversmooth textures, leading to the loss of critical image details.

To overcome the limitations of purely local and linear methods, early advances focused on transform-domain denoising techniques, focusing particularly on transforms that admit a sparse representation of the image. The wavelet transform provides a multi-resolution representation \cite{wavelet_mallat_book}, well-suited for modeling the piecewise smooth nature of natural images. The wavelet shrinkage framework \cite{donoho_soft_thresholding,donoho_johnstone} models noise attenuation by thresholding the wavelet coefficients. If \( W \) denotes an orthonormal wavelet transform and \( u = W z \) are the noisy coefficients, then a typical wavelet soft-thresholding scheme estimates the clean coefficients as $\hat{u}_i = \text{sign}(u_i) \cdot \max(|u_i| - \tau, 0)$, where \( \tau \) is a threshold chosen to balance noise removal and detail preservation. The denoised image is obtained by applying the inverse wavelet transform to the thresholded coefficients. Variants such as soft- and hard-thresholding \cite{donoho_soft_thresholding,donoho_johnstone}, as well as Bayesian thresholding rules \cite{wavelet_thresh_bayesian}, have been developed to adapt thresholds to local statistics, improving performance under varying noise levels.

Despite the success of transform-based methods, they struggle to fully exploit the inherent spatially repeating structures common in natural images. This limitation motivated the development of non-local algorithms that explicitly model self-similarity. A landmark example is the Non-Local Means (NLM) algorithm \cite{nlm_main}, which estimates each pixel as a weighted average of pixels across the entire image, with weights determined by the similarity between local neighborhoods. Formally, the estimate for a pixel at location \( i \) is given by $\hat{x}_i = \sum_{j} w_{ij} z_j$, where $w_{ij} = \frac{1}{S_i} \exp\left( -\frac{\| P_i - P_j \|_2^2}{h^2} \right)$, and \( P_i \) and \( P_j \) are patches centered at pixels \( i \) and \( j \) respectively. The parameter \( h \) controls the decay of the similarity function, and \( S_i \) is a normalizing constant ensuring that the weights sum to one. By exploiting repeated textures and patterns, NLM preserves fine structures that local or transform-based methods often miss.

Building upon the concept of non-local self-similarity, block-matching and collaborative filtering approaches emerged, with BM3D \cite{bm3d_main} becoming one of the most influential practical algorithms for image denoising. BM3D extends the non-local means principle by grouping similar patches into 3D stacks, applying collaborative transform-domain filtering within each group, and aggregating the estimates back to the image domain. The procedure involves block matching, a 3D linear transform (typically a combination of wavelet and discrete cosine transforms), hard-thresholding or Wiener filtering, and inverse transformation. This collaborative filtering step effectively separates signal and noise in the transform domain, and the aggregation of overlapping patches helps to reduce artifacts and improve robustness against mismatches in patch grouping.

In parallel, variational methods have provided a rigorous mathematical framework for denoising. Total Variation (TV) regularization \cite{RUDIN1992259_tv,Chambolle2004_tvmin} is a classic example formulated as $\hat{x} = \arg \min_{x} \frac{1}{2} \| y-x \|_2^2 + \lambda \| \nabla x \|_1$, where
\( \nabla x \) denotes the discrete image gradient and the \( \ell_1 \)-norm promotes sparsity in the gradient domain. TV denoising preserves edges by favoring piecewise constant regions while suppressing small oscillations due to noise. Despite its edge-preserving properties, TV regularization can suffer from the well-known \textit{staircasing} effect, where smooth intensity transitions are replaced by piecewise flat regions. To address this, higher-order regularization models and non-convex penalties have been proposed \cite{higher_order_tv,tgv,non_cvx_denoising}, providing more flexibility in capturing image textures and fine details.

As large datasets and increased computational resources became available, data-driven denoisers emerged as a dominant approach. Early learning-based methods focused on dictionary learning and sparse coding, where an overcomplete dictionary \( B \) is learned from (possibly noisy) image patches. Given a set of noisy patches \( z_i \), $1\leq i\leq p$, extracted from an image, their corresponding sparse codes \( \boldsymbol{\alpha}_i \) are estimated by solving
\[
\left(B,\left(\boldsymbol{\alpha}_i\right)_{i=1}^p\right) = \argmin_{B,\left(\boldsymbol{\alpha}_i\right)_{i=1}^p}\sum_{i=1}^p \left[\| z_i - B \boldsymbol{\alpha}_i \|_2^2 + \lambda_i \| \boldsymbol{\alpha}_i \|_1\right],
\]
and the denoised patches are reconstructed as \( B \boldsymbol{\alpha}_i \). The K-SVD algorithm \cite{ksvd_denoising,analysis_ksvd} is a notable example of this paradigm, providing a flexible and interpretable model that adapts to local structures.

With the advent of deep learning, convolutional neural networks (CNNs) have become the de facto standard for state-of-the-art image denoising. CNN-based models exploit large-scale training data to learn highly expressive and powerful mappings from noisy to clean images. A representative architecture, such as DnCNN \cite{dncnn_main}, employs multiple convolutional layers, batch normalization, and residual learning to directly estimate the noise component, which is then subtracted from the noisy input. The learned mapping can be described as $\hat{x} = z - f_{\theta}(z)$, where \( f_{\theta} \) denotes the denoising function parameterized by the learnable network weights \( \theta \). Residual learning accelerates convergence and stabilizes training by focusing the network on learning the noise distribution rather than the clean signal itself. Such models achieve remarkable generalization performance across a wide range of noise levels, structures, and image content.

Recent advances have further enhanced deep learning-based denoisers through the incorporation of attention mechanisms and transformer architectures \cite{swinir,restformer}, which capture long-range dependencies more effectively than purely local convolutions. These models can model non-local interactions within an image at a global scale, leading to improved reconstruction of repetitive patterns and structures. Moreover, self-supervised denoising approaches have gained significant traction, especially in applications where clean ground truth images are difficult or impossible to obtain. Methods such as Noise2Noise \cite{noise2noise} and Noise2Void \cite{noise2void} exploit the statistical independence of noise realizations or employ blind-spot training to learn denoisers directly from noisy data, broadening the practical applicability of learning-based denoising. Stein's unbiased risk estimation (SURE) approach \cite{sure_baraniuk} and equivariant denoising \cite{liu2025rotationequivariantselfsupervisedmethodimage} offer two attractive frameworks for self-supervised image denoising. While SURE replaces the mean squared error (MSE) with an unbiased estimate of it to eliminate the dependence on reference ground-truth images, the equivariance-based approach exploits rotational (or other) symmetries of images to achieve self-supervised learning of denoisers. 

The rise of generative models again provides new perspectives for image denoising. Denoising diffusion probabilistic models (DDPMs) \cite{NEURIPS2020_ddpm,improved_ddpm} represent one of the most promising directions. These models learn to reverse a forward diffusion process (represented through a stochastic differential equation (SDE)) that gradually adds noise to a clean image, enabling the generation of high-fidelity samples through a sequence of denoising steps. The iterative nature of diffusion models aligns naturally with the iterative refinement inherent in many classical denoising algorithms, making them a compelling candidate for plug-and-play priors in more complicated inverse problem settings. Although computationally intensive, diffusion-based denoisers have demonstrated state-of-the-art results and offer principled uncertainty quantification for the recovered image.

Throughout these decades of progress, one recurring theme has been the interplay between model-based image priors and data-driven learning of denoisers. Classical algorithms offer interpretability, well-defined mathematical properties, and provable convergence guarantees, but often lack the representational power needed to capture the complexity of natural images. Learned denoisers excel at modeling rich high-dimensional distributions but raise questions about stability, generalization, and robustness under distribution shifts. This tension has inspired the design of \rev{hybrid approaches, which embed} powerful learned denoisers into model-based optimization frameworks. This idea gave rise to the family of plug-and-play (PnP) methods, which have become highly successful because they bring together two advantages: the flexibility and performance of learned denoisers, and the interpretability and control offered by classical iterative schemes.

In summary, the evolution of image denoising algorithms reflects a remarkable trajectory, from simple linear filters to sophisticated non-local, transform-domain, variational, and deep learning-based methods. Each generation has expanded our understanding of natural image statistics and improved our ability to suppress noise while preserving important details in the image. These advances have laid the groundwork for modern inverse problem frameworks that leverage powerful denoising priors as modular components. As generative modeling and self-supervised learning continue to mature, they promise to inspire the next wave of innovation in imaging inverse problems and beyond.

\subsection{Inverse Problems and Regularization}
The study of regularization methods for ill-posed inverse problems in imaging originates from Hadamard's notion of well-posed problems, in the sense that the solution must exist, be unique, and vary continuously with respect to the observed data. The canonical linear inverse problem seeks to recover an image $x \in \mathcal{X}$ from noisy measurements $y \in \mathcal{Y}$ related through the forward model:
\begin{equation}
y = Kx + w,
\end{equation}
where $K: \mathcal{X} \rightarrow \mathcal{Y}$ is a compact operator between Hilbert spaces, and $w$ represents measurement noise. The ill-posed nature manifests through three distinct aspects of the operator-theoretic framework: either \textit{injectivity} or \textit{surjectivity} of the forward operator may not hold, or \textit{stability} of the solution map might be violated. For instance, if $K$ is a compact operator with an infinite-dimensional range, then surjectivity and stability are violated. This is, for instance, the case for the ray transform operator that underlies many medical imaging modalities, such as computed tomography (CT) and positron emission tomography (PET). To see this, first consider the case where the range $\mathscr{R}(K)$ of the forward operator $K$ is not closed, implying that solutions may not exist for an arbitrary $y \in \mathcal{Y}$, as any non-zero noise component $w$ orthogonal to $\mathscr{R}(K)$ renders the problem unsolvable in the strict sense. Second, the potential non-triviality of $\mathscr{N}(K)$, the null-space of the forward operator $K$, violates uniqueness, particularly evident in limited-angle tomography where certain features become invisible. Most critically, the unboundedness of the generalized inverse $\rev{K}^\dagger$, when restricted to $\mathscr{R}(K)$, leads to extreme sensitivity to noise, leading to an unstable solution (i.e., a small amount of noise in the measurement yields a drastically different solution).

\noindent These challenges become explicit through the singular value decomposition of the compact operator $K = \sum_{m=1}^\infty \sigma_m \langle \cdot, u_m \rangle v_m$, where the asymptotically vanishing singular values $\sigma_m \rightarrow 0$ cause the naive solution, given by
\begin{equation}
x^\dagger = \sum_{m=1}^\infty \frac{\langle y, v_m \rangle}{\sigma_m} u_m,
\end{equation}
cause noise components corresponding to small singular values to be catastrophically amplified. Regularization theory addresses this instability by constructing families of approximate solutions ${x_\alpha}$ through variational formulations:
\begin{equation}
x_\alpha = \argmin_{x \in \mathcal{X}} \frac{1}{2}\|Kx - y\|_\mathcal{Y}^2 + \alpha R(x),
\label{eq_varreg_basic}
\end{equation}
where the functional $R: \mathcal{X} \rightarrow \mathbb{R} \cup \{+\infty\}$ encodes prior knowledge about plausible solutions. \rev{A solution $x_{\alpha}$ to \eqref{eq_varreg_basic} is alternatively denoted as the output of a parametric reconstruction operator $\mathcal{R}_\alpha \coloneqq \mathcal{R}(\cdot,\alpha):\mathcal{Y}\rightarrow \mathcal{X}$.} The classical Tikhonov regularization \cite{tikhonov_main} employs $R(x) = \rev{\frac{1}{2}}\|x\|_\mathcal{X}^2$, yielding the explicit solution $x_\alpha = (K^\top K + \alpha I)^{-1}K^\top y$, \rev{where $K^\top$ denotes the adjoint\footnote{{We use the notation $K^\top$ to denote the adjoint of $K$ regardless of the image domain $\mathcal{X}$. When $\mathcal{X}=\mathbb{R}^n$, the adjoint operator $K^\top$ reduces to the transpose of the matrix $K$.}} of the forward operator $K$}. This solution corresponds to a spectral filter with coefficients $f_\alpha(\sigma) = \sigma/(\sigma^2 + \alpha)$. While the quadratic ($L^2$) penalty on $x$ guarantees stability of the reconstruction, it does so by uniformly penalizing deviations across all frequencies. This effect can be seen from the Euler–Lagrange optimality condition $K^\top(Kx - y) + \alpha x = 0$, where the regularization parameter $\alpha$ acts as a frequency-independent damping term. The outcome is a systematic suppression of both noise and fine-scale features, resulting in overly smoothed reconstructions and noticeable loss of sharp edge structures. To mitigate this limitation, one may instead consider sparsity-promoting penalties such as the $L^1$ norm \cite{lasso_tibshirani}, \rev{typically applied on a transform domain. For instance, to promote sparsity of the wavelet coefficients of the image, one can choose $R(x)=\|Wx\|_1$, where $W$ denotes an appropriate wavelet transform operator. Such sparsity-promoting $L^1$ norm-based regularizers are non-differentiable, leading to a non-smooth variational optimization problem for reconstruction (which requires iterative solvers, unlike the variational problem with the Tikhonov regularizer admitting a closed-form solution). Another popular and widely adopted choice is the TV regularizer $R(x)=\|\nabla x\|_1$ discussed in Section \ref{sec_brief_denoising_survey}, which effectively allows the reconstruction to tolerate large local variations and thus preserve sharp discontinuities such as edges. However, the piecewise-constant bias induced by TV regularization tends to replace smooth gradients in the original image by artificial flat regions separated by sharp transitions of intensity \cite{RUDIN1992259_tv,Chambolle1997_tv_related}.  Moreover, textures and fine oscillatory details tend to be lost, since they are not easily represented in the sparse model enforced on the gradient image by the TV regularizer.} 

Modern approaches combine multiple regularization functionals through formulations of the form $R(x) = \sum_{l=1}^L \alpha_l R_l(x)$, where typical components include sparsity-promoting terms $\|\Psi x\|_1$ in learned dictionaries, higher-order derivatives $\|\nabla^2 x\|_1$, and nonlocal operators capturing long-range image dependencies.

One of the most significant theoretical advancements in recent years has been the development of learned regularization through the plug-and-play framework \cite{venkat_pnp_6737048}, where advanced denoising operators $D_\sigma$ are interpreted as proximal operators corresponding to implicit regularizers:
\begin{equation}
D_\sigma(x) \approx \prox_{\sigma^2 R}(x) = \argmin_z \frac{1}{2}\|z - x\|^2 + \sigma^2 R(z).
\end{equation}
This interpretation leads to provably convergent algorithms (at least in the sense of fixed-point convergence) when the denoiser satisfies appropriate non-expansiveness conditions. The resulting methods combine the theoretical foundations of variational regularization with the excellent empirical performance of data-driven image priors induced by denoisers.

Current theoretical challenges include the rigorous characterization of the implicit regularizers associated with modern denoising architectures, the extension to nonlinear forward models, and the development of convergence rates under weaker assumptions on the denoising operators. These questions represent active areas of research at the intersection of functional analysis, optimization theory, and statistical learning. 

\subsection{Convergent Regularization}
To obtain stable solutions to inverse problems, a mechanism is needed to handle varying noise levels in the measurement. \rev{When the measurement noise level is large (small), one must apply a stronger (weaker) regularization: this ensures that the variational framework \eqref{eq_varreg_basic} for reconstruction optimally trades off data-fidelity with the prior knowledge through the parameter $\alpha$. The explicit dependence of $\alpha$ on the measurement noise can be explained by interpreting \eqref{eq_varreg_basic} as a Bayesian maximum a-posteriori estimation problem with an image prior proportional to $\exp(-\beta\,R(x))$ and Gaussian measurement noise with variance $\gamma_w^2$, resulting in $\alpha=\beta\gamma_w^2$.} For this purpose, the concept of convergent regularization has proven highly useful. Regularization can be roughly understood as a convergence requirement to a unique solution, such as the minimum-norm solution $x^\dagger$, where convergence occurs as the noise level $\delta \rightarrow 0$. Formally, consider the previously discussed reconstruction operator $\mathcal{R}_\alpha \coloneqq \mathcal{R}(\cdot,\alpha)$, which parametrizes a family of continuous operators $\mathcal{R}_\alpha:\mathcal{Y}\to \mathcal{X}$. The parameter $\alpha$ depends on the noise level $\delta>0$, where $\| y^{\delta} - y^0 \| \leq \delta$ and $y^0 \coloneqq K x$ denotes noise-free measurement data. We say that the family of reconstruction operators is a \textit{convergent regularization} method if there exists a parameter choice rule $\alpha = \alpha(\delta, y^{\delta})$ such that reconstructions $x^\delta \coloneqq \mathcal{R}_{\alpha(\delta,y^\delta)}(y^\delta)$ converge to $x^\dagger \coloneqq \rev{K}^{\dagger} y^0$ (given by the pseudo-inverse) as noise vanishes, in the sense that
\begin{equation}
	\limsup_{\delta \to 0} \left\| x^\delta
	-
	x^{\dagger}
	\right\|_{\mathcal{X}}=0
	\quad\text{as}\quad
	\limsup_{\delta \to 0} \,
	\{\alpha(\delta,y^{\delta})\}
	= 0.
    \label{eqn:classicRegCondition}
\end{equation}
In other words, we have point-wise convergence of the reconstruction operators to the pseudo-inverse, i.e. $\mathcal{R}_{\alpha(\delta, y^{\delta})}(y^\delta) \to \rev{K}^\dagger y^0$ as $\delta\to 0$. We refer interested readers to \cite{engl1996regularization,scherzer2009variational,benning2018modern} for a detailed discussion on convergent regularization schemes (and several convergence rate results in the classical regularization literature). While this is somewhat restrictive as it only considers convergence to the least-squares minimum-norm solution, this nevertheless serves as an important tool to design learned regularization methods, i.e., learned reconstruction approaches that formally satisfy the above convergence criterion. 

\section{Proximal Splitting Algorithms}\label{sec1}
The heart of PnP lies within monotone operator theory, particularly operator splitting. Informally, splitting methods solve a composite optimization problem using simpler gradient-like operations, interpreted as data fidelity steps and regularization steps within the PnP framework. We will introduce in this section the notion of proximal operators and their centrality in convex analysis, and demonstrate how convergence results in monotone operator theory relate to composite convex optimization and further to convergence of PnP methods. For a more in-depth exposition on convex analysis and monotone operator theory, we refer interested readers to \cite{bauschke2011convex}.

\rev{The origins of splitting methods can be traced back to the seminal work of Douglas and Rachford (1956) on solving heat conduction problems \cite{douglas1956numerical}. Lions and Mercier (1979) later generalized these ideas to maximal monotone operators in Hilbert spaces \cite{lions1979splitting}. The modern ADMM framework emerged through the work of Gabay, Mercier (1976) \cite{gabay1976dual}, and Glowinski (1985) \cite{glowinski1985numerical}, with Eckstein (1989) establishing the definitive connection to DRS \cite{eckstein1989splitting}. Recent advances have focused on several key directions, for instance, momentum-based variants incorporating Nesterov-type acceleration, stochastic implementations for large-scale problems, nonconvex extensions with convergence guarantees, and distributed implementations for multi-agent systems. The theoretical understanding of these methods continues to deepen, with new connections to differential inclusions and variational inequalities being actively explored. }

\subsection{Foundations of Proximal Calculus} 
The proximal operator, introduced by Moreau in 1962, serves as the cornerstone of modern nonsmooth optimization. Given a proper, closed, and convex function $f:\mathbb{R}^n \to \mathbb{R}\cup\{+\infty\}$, its proximal operator is defined through the solution of the following variational problem:
\begin{equation}\label{eq:proxDef}
{\prox}_{\lambda f}(v) = \argmin_x \left( f(x) + \frac{1}{2\lambda}\|x - v\|_2^2 \right).
\end{equation}
In the case where $f$ is differentiable, the proximal operator may be interpreted as a backwards Euler discretization of gradient flow $\dot{x}(t) = -\nabla f (x(t))$. With step size $\eta=\lambda$, the \rev{proximal scheme $x_{k+1} = \prox_{\lambda f}(x_k)$ is equivalently} given by $x_{k+1} = x_k - \eta \nabla f(x_{k+1})$, which \rev{may be shown using} the first-order optimality condition \rev{in} \labelcref{eq:proxDef}. The proximal operator generalizes this concept to nonsmooth and $\infty$-valued functions. Furthermore, the proximal operator is a generalization of projection onto a convex set. Suppose $f = \chi_C$ is the characteristic function of a convex set $C$, taking values $0$ in $C$ and $+\infty$ otherwise. Then, the proximal operator $\prox_{\lambda f}$ is precisely the (Euclidean) projection onto $C$. This equivalence is useful when interpreting constrained optimization.

\noindent For a class of general convex functions, the proximal operator has several useful functional properties, as stated in the following proposition.
\begin{proposition}[{\cite{rockafellar1997convex, ekeland_book, rockafellar2009variational}}]\label{prelim:prop:proxProps}
    For a proper closed convex function $f$, the proximal operator is well-defined and is single-valued. Moreover, it satisfies the following:
    \begin{enumerate}
        \item $\prox_f$ is non-expansive (i.e., 1-Lipschitz) and continuous. 
        \item Fixed points of $\prox_f$ correspond to minimizers of $f$: 
            \[\{x_0 \in {\R^n}\mid x_0 = \prox_f(x_0)\} = \argmin_{x\in \R^n} f(x).\]
        \item (Moreau's identity) $\prox_f + \prox_{f^*} = \id$, where $\id$ is the identity map on $\R^n$. 
    \end{enumerate}
\end{proposition}
\subsubsection{Monotone Operators}
In addition to its functional properties, the proximal operator's importance to imaging stems from its deep connection to the calculus of subgradients, which can be understood through the lens of monotone operators. In what follows, we give a brief review of monotone operators, paving the way for discussing splitting algorithms that underlie the modern PnP methods.  
\begin{definition}[Monotonicity]
    A set-valued mapping $T:\R^n \rightrightarrows \R^n$ is said to be $\emph{monotone}$ if for all $x, x' \in \R^n,\, p \in T(x),\, p' \in T(x')$,
    \[\langle p-p', x-x' \rangle \ge 0,\]
    and \emph{strictly monotone} if the inequality is strict for $x \ne x'$. The \emph{resolvent} of $T$ is the operator $J_T \coloneqq (\id + T)^{-1}$, and the \emph{reflected resolvent} is $R_T \coloneqq 2J_T - \id$.
    A set-valued mapping $T$ is said to be \emph{maximally monotone} if its graph $\Gamma(T) := \{(x, p): x \in \R^n,\, p \in T(x)\}$ is not contained within the graph of another monotone operator.
\end{definition}
\noindent For a (proper and closed) convex function $f$, the subdifferential operator $\partial f$ can be seen to be a monotone operator (and indeed, something stronger called \rev{maximally} cyclically monotone) \cite{rockafellar1966characterization}. Analogously to the backward Euler interpretation above, the proximal map can be equivalently characterized as the resolvent of the subdifferential operator \cite[Sec. 12.C.]{rockafellar2009variational}:
\begin{equation}
\prox_{\lambda f} = J_{\rev{\lambda}\partial f} \coloneqq (\id + \lambda\, \partial f)^{-1}.
\end{equation}
%
\rev{By Minty's theorem, we have that the resolvent of a monotone operator is defined everywhere if and only if the monotone operator is maximally monotone \cite[Thm. 21.1]{bauschke2011convex}. Moreover, a function $T:\R^n \rightarrow \R^n$ is firmly non-expansive, i.e. $\|T(x-y)\|^2 + \|(\id-T)(x-y)\|^2 \le \|x-y\|^2$, if and only if it is the resolvent of a maximally monotone operator \cite[Cor. 23.8]{bauschke2011convex}. This firm non-expansiveness condition is sufficient for the fixed point iteration $x_{n+1} = Tx_n$ to converge.

Other interesting formulas can be reformulated in terms of maximal monotone operators. For example, Moreau's identity $\prox_f + \prox_{f^*} = \id$ can be reformulated using resolvents as $\id = J_{\gamma A} + \gamma^{-1} J_{\gamma^{-1}A^{-1}} \circ \gamma^{-1} \id$ \cite[Prop. 23.18]{bauschke2011convex}. Another lesser known identity states: if $f$ is proper closed and convex, and ${}^\gamma\! f \coloneqq f \circ \prox_f$ is its Moreau envelope, then $\prox_{{}^\gamma\! f} (x) = x + (\gamma+1)^{-1} (\prox_{(\gamma+1 )f}(x) - x)$ \cite[Prop. 23.29]{bauschke2011convex}.} These theoretical foundations \rev{and relations} underpin the development of proximal algorithms, particularly the proximal gradient method, arising in convex optimization by leveraging results from monotone operator theory.


\subsubsection{Composite Optimization and Operator Splitting}
Common optimization problems encountered in variational image recovery take the composite form $\argmin_x f(x) + g(x)$, where $f, g$ are proper closed convex functions with possibly different regularity conditions. Considering the first-order optimality conditions, the composite optimization problem is equivalent to the monotone inclusion problem $0 \in Ax + Bx$, where $A=\partial f$ and $B=\partial g$ are both maximally monotone operators, arising naturally while finding a minimizer of the sum of two convex functions. 

For general monotone operators $A$ and $B$, one possible approach is to consider root solving using the resolvent $J_{\lambda(A+B)}$. However, this may be difficult to compute. For example, taking $A=\partial f$ and $B=\partial g$, the subgradients of $f$ and $g$ respectively, this is equivalent to computing the proximal operator of the composite function $f+g$. Therefore, one seeks to find a zero of $A+B$, using only their resolvents $J_{\lambda A}$ and $J_{\lambda B}$. This is useful in the context of convex problems where $f$ and/or $g$ have easily computable proximals, while $f+g$ does not. This process of splitting the resolvent of $A+B$ into the resolvents of its components is referred to as a \textit{splitting algorithm} and can be done in different ways \cite{lions1979splitting}. We present two simple versions, which are by far the most widely employed splitting techniques in convex optimization: \textit{proximal gradient descent} (sometimes referred to as \textit{forward-backward splitting} (FBS)) and the \textit{Douglas--Rachford splitting} (DRS) \cite{douglas1956numerical}. 

\rev{In the following results, we denote the fixed points of an operator $A$ by $\Fix A \coloneqq\{x \in \mathcal{X} \mid Ax = x\}$, and the zeros of a (possibly set-valued) mapping $B$ by $\zer B \coloneqq \{x \in \mathcal{X} \mid 0 \in Bx\}$.}


\begin{theorem}[{Forward-backward algorithm \cite[Cor. 27.9]{bauschke2011convex}}]\label{thm:FBS}
    For a Hilbert space $\mathcal{X}$, let $A:\mathcal{X} \rightarrow \mathcal{X}$ be $\beta$-cocoercive (i.e. $\langle x-y, Ax-Ay\rangle \ge \beta \|Ax-Ay\|^2$) for some $\beta>0$, and $B:\mathcal{X} \rightrightarrows \mathcal{X}$ be maximally monotone. Let $\lambda \in (0,2\beta)$ \rev{and $x_0 \in \R^n$ be an initialization}, and further set $\delta = \min\{1,\,\beta/\lambda\}+1/2$. Suppose that $\zer(A+B)\ne \emptyset$, and define the forward-backward iterations as follows, 
    \begin{equation}\label{eq:FBS}
        \begin{cases}
            y_k = x_k - \lambda Ax_k, \\
            x_{k+1} = J_{\lambda B}y_k.
        \end{cases}
    \end{equation}
    The iterates satisfy the following:
    \begin{enumerate}
        \item $(x_k)_{k \in \mathbb{N}}$ converges weakly to a point in $\zer(A+B)$;
        \item Suppose $x \in \zer(A+B)$. Then $Ax_n$ converges strongly to $Ax$.
    \end{enumerate}
\end{theorem}

\begin{theorem}[{Douglas--Rachford Splitting \cite[Thm. 25.6]{bauschke2011convex}}]\label{thm:DRS}
    For a Hilbert space $\mathcal{X}$, let $A, B:\mathcal{X} \rightrightarrows \mathcal{X}$ be maximally monotone operators such that $\zer(A+B) \ne \emptyset$. Let $\lambda>0$ be a step size and $x_0 \in \rev{\mathcal{X}}$ be an initialization. Consider the iterations \rev{
    \begin{equation}\label{eq:DRS}
        \begin{cases}
            y_k = J_{\lambda A} x_k,\\
            z_k = J_{\lambda B} (2y_k - x_k),\\
            x_{k+1} = x_k + z_k-y_k;\\
        \end{cases}
    \end{equation}}
    which can be expressed more succinctly as,\rev{
    \begin{equation}\label{eq:DRSUpdate}
        x_{k+1} = J_{\lambda B}(2 J_{\lambda A} - \id) x_n+ (\id - J_{\lambda A}) x_k.
    \end{equation}}
    Then, there exists a fixed point $x \in \Fix R_{\lambda B}R_{\lambda A}$ such that the following hold:
    \begin{enumerate}
        \item $J_{\lambda A}(x) \in \zer(A+B)$,
        \item $y-z_k$ converges strongly to zero,
        \item $x_k$ converges weakly to $x$, and
        \item $y$ and $z_k$ converge weakly to $J_{\lambda A}(x)$.
    \end{enumerate}
\end{theorem}
\noindent Note that in the case where the Hilbert space $\mathcal{X}$ is finite-dimensional, weak convergence is equivalent to strong convergence. \rev{The reflected resolvent of $A$ appears in \labelcref{eq:DRSUpdate}; since $A$ is maximally monotone, the reflected resolvent is a non-expansive operator \cite[Cor. 23.10]{bauschke2011convex}, allowing for a contraction-like argument.}

To obtain the corresponding \rev{optimization method}, simply let $A$ and $B$ be the subdifferentials of some proper closed convex functions $f$ and $g$. We get convergence to a zero of \rev{${A+B}$}, equivalently a fixed point of $\prox_{f+g}$, using only proximal operators or subgradients of $f$ and $g$ separately. Furthermore, the fixed point is a minimum of $f+g$. 

\rev{FBS and DRS impose different requirements for efficient computation, where the former needs that $\nabla f$ and $\prox_g$ are easy to compute (and $\nabla f$ is Lipschitz, hence cocoercive by Baillon--Haddad), and the latter requires that both $\prox_f$ and $\prox_g$ are easy to compute. Moreover, DRS can be extended to finding a zero of a finite sum of maximally monotone operators, with the resulting algorithm known as the \emph{parallel splitting algorithm} \cite[Prop. 25.7]{bauschke2011convex}. Sharp convergence rates for FBS and DRS applied to composite optimization can be found in \cite{davis2015convergence,molinari2019convergence}.}



\subsubsection{PnP Proximal Gradient Descent}
\rev{By casting the above monotone inclusion problem \labelcref{eq:FBS} in the scope of convex functions, with $A$ being a derivative and $B$ being a proximal operator, we can obtain splitting schemes that optimize the sum of two convex functions, where one of the functions is smooth. Letting $A=\nabla f$, $B=\partial g$}, we consider \rev{solving the following} composite optimization problem:
\begin{equation}
\min_x f(x) + g(x)
\end{equation}
where $f$ is $L$-smooth and $g$ admits efficient proximal evaluations. Note that $L$-smoothness of $f$ corresponds to $1/L$-cocoercivity of $\nabla f$ using the Baillon--Haddad theorem \cite{baillon1977quelques}. The \textit{proximal gradient descent} (PGD) method generates iterates via:
\begin{equation}\label{eq:PGD}
x_{k+1} = \prox_{\lambda g}\left(x_k - \lambda \nabla f(x_k)\right)
\end{equation}
\noindent The convergence properties of this scheme are characterized below.
\begin{theorem}[Proximal Gradient Convergence {\cite[Sec. 10]{beck2017first}}]
Let $f$ be $\mu$-strongly convex and $L$-smooth for some $\mu\ge 0, L>0$, and let $g$ be convex. For step size $\lambda \in (0, 2/L)$ \rev{and initialization $x_0$}, the PGD iterates satisfy
\begin{equation}
\|x_k - x^*\|_2 \leq \varrho^k \|\rev{x_0} - x^*\|_2,
\end{equation}
where $\varrho = \max(|1 - \lambda L|, |1 - \lambda \mu|) < 1$, \rev{and $x^*$ is the minimizer of $f+g$}. When $f$ is merely convex, the objective error decays as $\mathcal{O}(1/k)$. Moreover, the minimum of the residuals satisfies
\begin{equation}
    \min_{l \le k} \|x_l - x_{l+1}\| = \mathcal{O}(1/k).
\end{equation}
\end{theorem}
\noindent The PnP-PGD method is obtained by replacing the $\prox_{\lambda g}$ term in \labelcref{eq:PGD} with a (Gaussian) denoiser $D_{\sigma}$, where $\sigma$ denotes the  standard deviation of noise that the denoiser can eliminate: 
\begin{equation}\label{eq:PnP-PGD}\tag{PnP-PGD}
    x_{k+1} = D_\sigma\left(x_k - \lambda \nabla f(x_k)\right).
\end{equation}
\subsubsection{Relaxed PnP Proximal Gradient Descent}
The PGD iterations can be relaxed to accommodate for weakly convex functions $g$ while minimizing $f+g$. The relaxed iterations, for \rev{some relaxation parameter} $\alpha\in (0,1)$, are \cite{hurault2024convergent,chambolle2016ergodic}:
\begin{equation} \label{eqs:aPGD}
    \begin{cases}
        q_{k+1} = (1-\alpha)x_k + \alpha y_k,\\
        y_{k+1} = \prox_{\lambda g} (y_k - \lambda \nabla f(q_{k+1})),\\
        x_{k+1} = (1-\alpha) x_k + \alpha y_{k+1}.
    \end{cases}
\end{equation}
The relaxed PGD algorithm, also known as $\alpha$PGD, enjoys similar convergence results as FBS and DRS. Moreover, the convergence theory can handle weakly convex functions, albeit with a modified objective functional involving an additional residual term. \rev{}

\begin{theorem}[{\cite[Thm. 2]{hurault2024convergent}}]
    Let $f$ be convex and \rev{$L_f$}-smooth, and $g$ be $M$-weakly convex. Then for $\alpha \in (0,1)$ and $\lambda < \min((\alpha L_f)^{-1}, \alpha M^{-1})$, define the objective \rev{$F =  f + g$}. The iterates $x_k$ in \labelcref{eqs:aPGD} satisfy
    \begin{enumerate}
        \item $F(x_k) + \frac{\alpha}{\rev{2}} \left(1- \frac{1}{\alpha}\right)^2 \|x_k - \rev{x_{k-1}}\|^2$ is non-increasing and convergent;
        \item The sequence $(x_k)$ has finite length, \rev{i.e., $\sum_{k \in \mathbb{N}} \|x_{k+1} - x_k\| < +\infty$. Moreover, }$\min_{l<\rev{k}}\|x_{l+1}-x_l\| = \mathcal{O}(1/\sqrt{k})$;
        \item The cluster points of the sequence $(x_k)$ are stationary points of $F$.
    \end{enumerate}
\end{theorem}
\noindent Notably, the gap between the minimum residual rate of $\mathcal{O}(1/\sqrt{k})$ for $\alpha$PGD and $\mathcal{O}(1/k)$ for PGD is a consequence of convexity. In PnP methods, we usually deal with weakly convex functions, for which the $\mathcal{O}(1/\sqrt{k})$ rate comes from \rev{the smooth convex fidelity term}. Analogous to PnP-PGD, the relaxed PnP-$\alpha$PGD method arises from replacing the proximal with a denoiser in \labelcref{eqs:aPGD},
\begin{equation} \label{eqs:PnP_aPGD} \tag{PnP-$\alpha$PGD}
    \begin{cases}
        q_{k+1} = (1-\alpha)x_k + \alpha y_k,\\
        y_{k+1} = D_\sigma (y_k - \lambda \nabla f(q_{k+1})),\\
        x_{k+1} = (1-\alpha) x_k + \alpha y_{k+1}.
    \end{cases}
\end{equation}
\subsubsection{PnP Douglas--Rachford Splitting}
For proper, convex, and closed functions $f$ and $g$, one can substitute $\rev{A=\partial f}, B=\partial g$ into \labelcref{eq:DRS} to yield the update
\begin{equation}\label{eq:DRSiter}
    x_{k+1} = \prox_{\lambda g}(2 \prox_{\lambda f} - \id) x_k + (\id - \prox_{\lambda f}) x_k.
\end{equation}

\begin{theorem}[{\cite[Thm. 3.1]{he2015convergence}}]
    The residuals in the DRS iteration \rev{\labelcref{eq:DRSiter}}, given by 
    \begin{equation*}
        e_k = x_{k} - \prox_{\lambda g}(x_k - \lambda \nabla f(x_k))
    \end{equation*}
    decay as $\|e_k\|^2 = \mathcal{O}(1/k)$.
\end{theorem}
\noindent The asymmetry of the Douglas--Rachford splitting gives rise to two possible splittings by switching the roles of $f$ and $g$ \cite{hurault2022proximal}. \rev{These two splittings} give rise to PnP-DRS and PnP-DRSdiff, defined as follows. Notably, the two variants of PnP with the DRS have slightly different convergence assumptions on $g$, \rev{and with PnP-DRSdiff further requiring that $f$ is differentiable}. The PnP-DRS iterations are given by
\begin{equation} \label{eqs:PnP_DRS} \tag{PnP-DRS}
    \begin{cases}
        y_{k+1} = D_\sigma(x_k),\\
        z_{k+1} = \prox_{\lambda f}(2y_{k+1} - x_k),\\
        x_{k+1} = x_k + (z_{k+1} - y_{k+1});
    \end{cases}
\end{equation}
while the PnP-DRSdiff algorithm generates the following iterates:
\begin{equation} \label{eqs:PnP_DRSdiff}  \tag{PnP-DRSdiff}
    \begin{cases}
        y_{k+1} = \prox_{\lambda f}(x_k),\\
        z_{k+1} = D_\sigma(2y_{k+1} - x_k),\\
        x_{k+1} = x_k + (z_{k+1} - y_{k+1}).
    \end{cases}
\end{equation}
\subsection{ADMM: Constrained Optimization to Operator Splitting}
While PGD and DRS consider composite optimization in one variable, a more general problem may include equality constraints, such as in a Lagrangian. Alternatively, one may be interested in a problem of the form $f(x) + g(Kx)$ for some forward operator $K$ with a computable adjoint $K^\top$. The \textit{alternating direction method of multipliers} (ADMM) algorithm \cite{boyd_admm_book} is equipped to address problems with a separable structure having the \rev{more} general form 
\begin{equation}
\min_{x,z} f(x) + g(z) \quad \text{subject to} \quad Kx + K'z = c,
\end{equation}
\rev{where $K'$ is another linear operator, and the slack variable $c$ may arise from e.g. data constraints arising from an underlying inverse problem.} ADMM can be derived from applying DRS to a dual formulation \cite{gabay1983chapter,ryu2016primer}. To see this, introduce the augmented Lagrangian according to the equality constraint:
\begin{equation}
\mathcal{L}_\rho(x,z,u) = f(x) + g(z) + u^\top(Kx + K'z - c) + \frac{\rho}{2}\|Kx + K'z - c\|_2^2.
\end{equation}
ADMM then generates iterates through alternating minimization:
\begin{align}
\begin{cases}
    x_{k+1} = \argmin_x \, \mathcal{L}_\rho(x,z_k,u_k), \\
z_{k+1} = \argmin_z \, \mathcal{L}_\rho(x_{k+1},z,u_k), \\
u_{k+1} = u_k + \rho(Kx_{k+1} + K'z_{k+1} - c).
\end{cases}
\label{eq_admm_iter}
\end{align}
The connection to Douglas--Rachford splitting emerges when considering the dual problem. We demonstrate this in the special case where $K' = -\id$ and $c = 0$. In this case, the primal problem becomes
\begin{equation}\label{eq:ADMMSplitprimal}
\min_{x,z} f(x) + g(z) \quad \text{subject to} \quad Kx = z,
\end{equation}
and the corresponding dual problem is given by
\begin{equation}
\max_u -f^*(-K^\top u) - g^*(u),
\end{equation}
where $K^\top$ is the adjoint of the operator $K$, equivalently, the matrix transpose when \rev{$\mathcal{X}=\mathbb{R}^n$}. The optimality conditions for the dual problem are:
\begin{align*}
0 &\in -\rev{K}\partial f^*(-K^\top u^*) + \partial g^*(u^*) \\
&\Updownarrow \notag\\
\exists x^*, z^* \ \text{s.t.}\quad  z^* &= Kx^* \quad \text{where} \quad x^* \in \partial f^*(-\rev{K^\top} u^*), \, z^* \in \partial g^*(u^*).
\end{align*}
Define the maximally monotone operators:
\begin{align*}
T_1(u) &= -K\partial f^*(-K^\top u) = \partial(f^* \circ (-K^\top))(u), \\
T_2(u) &= \partial g^*(u).
\end{align*}
We need to solve the inclusion problem $0 \in T_1(u) + T_2(u)$. Using DRS \rev{\labelcref{eq:DRS} on the primal variable $z_k+w_k$, and where the dual variable is given by $u_k$}, the iterations may be rewritten as 
\begin{align*}
    u_{k+1} &= J_{\rho T_1}(z_k + w_k),\\
    z_{k+1} &= J_{\rho T_2} (u_{k+1} - w_k),\\
    w_{k+1} &= w_k + z_{k+1} - u_{k+1}.
\end{align*}
The resolvent of $T_2$ may be recognized exactly as a proximal operator $J_{\rho T_2} = J_{\rho \partial g^*} = \prox_{\rho g^*} =\id - \rho \prox_{\rho^{-1} g}(\rho^{-1} \cdot)$. The resolvent of the first operator may also be computed as 
\begin{align*}
    J_{\rho T_1}(v) &= \prox_{\rho f^* \circ (-K^\top)}(v) \\
    &= v - \rho \prox_{\rho^{-1} [f^* \circ (-K^\top)]^*}(\rho^{-1}v)\\
    &= v - \rho \argmin_x \inf_{z \rev{\text{ s.t.}} -Az = x} f(z) + \frac{\rho}{2} \|x - \rho^{-1}v\|^2\\
    &= v + \rho K \argmin_z\left( f(z) + \frac{\rho}{2}\|\rev{K}z + \rho^{-1} v\|^2\right).
\end{align*}
Therefore, we have that:
\begin{align*}
    u^+ &= z+w+ \rho K\hat{x}_1,\quad \text{where}\\
    \hat{x}_1 &= \argmin_{x_1}\left(f(x_1) + \frac{\rho}{2} \|Kx_1 + \rho^{-1}(z+w)\|^2 \right)\\
    &= \argmin_{x_1}\left(f(x_1) + z^\top (Kx_1)+\frac{\rho}{2} \|Kx_1 + \rho^{-1}w\|^2 \right),
\end{align*}
and 
\begin{align*}
    z^+ &= \prox_{\rho g^*} (u^+ - w) = \prox_{\rho g^*} (z + \rho K \hat{x}_1)\\
    &= z + \rho K \hat{x}_1 - \rho \prox_{\rho^{-1} g} (\rho^{-1} (z+\rho K\hat{x}_1))\\
    &= z + \rho(K \hat{x}_1 - \hat{x}_2), \quad \text{where}\\
    \hat{x}_2 &= \prox_{\rho^{-1} g} (\rho^{-1} (z+\rho K\hat{x}_1))\\
    &= \argmin_{x_2} \left(g(x_2) + \frac{\rho}{2} \|x_2 - K \hat{x}_1- \rho^{-1}z\|^2\right)\\
    &= \argmin_{x_2} \left(g(x_2) - z^\top x_2+ \frac{\rho}{2} \|K \hat{x}_1- x_2\|^2\right).
\end{align*}
Finally, the dual DRS update $w^+ = w + z^+ - u^+$ simplifies to $w^+ = -\rho \hat{x}_2$. This may be re-substituted into the expression for $\hat{x}_1$. Now considering the variables $(\hat{x}_{k,1}, \hat{x}_{k,2}, z_k)$, the iteration simplifies to
\begin{equation}\label{eqs:ADMMhats}
    \begin{cases}
        \hat{x}_{k+1, 1} =  \argmin_{x_1}\left(f(x_1) + z_k^\top (Kx_1)+\frac{\rho}{2} \|K\rev{x_{1}} - \hat{x}_{k,2}\|^2 \right),\\
        \hat{x}_{k+1, 2} =  \argmin_{x_2}\left(g(x_2) - z_k^\top x_2+\frac{\rho}{2} \|K\rev{\hat{x}_{k+1,1} - x_2}\|^2 \right),\\
        z_{k+1} = z_k + \rho(K \hat{x}_{k+1,1} - \hat{x}_{k+1,2}).
    \end{cases}
\end{equation}
Recalling the (simplified) form of the augmented Lagrangian
\begin{equation}\label{eq:ADMMAugLagrang}
\mathcal{L}_\rho(x_1,x_2,z) = f(x_1) + g(x_2) + z^\top(Kx_1 - x_2) + \frac{\rho}{2}\|Kx_1 - x_2 \|_2^2,
\end{equation}
We observe that the \rev{minimization steps in \labelcref{eqs:ADMMhats} are precisely the minimizations with respect to the first two arguments of \labelcref{eq:ADMMAugLagrang}}. \rev{Such an equivalence provides a powerful link to the theoretically simpler DRS, providing initial convergence results for ADMM with relaxation factors $\rho$ through the monotone operator framework \cite{eckstein1992douglas}.} In the context of imaging inverse problems, the PnP variant of the ADMM algorithm emerges by replacing the proximal operator (with respect to a nonsmooth regularizer $g$) with an off-the-shelf denoiser $D$, while minimizing the variational objective with an $\ell_2^2$ fidelity term: $\underset{x}{\min}\,\frac{1}{2}\|y-Kx\|_2^2+g(x)$. To solve this problem using ADMM, one applies the variable separation trick to reformulate the problem as 
\begin{equation}
    \underset{x,z}{\min}\,\frac{1}{2}\|y-Kx\|_2^2+g(z) \text{\,\,subject to\,\,}x=z. 
    \label{eq_var_opt_for_admm}
\end{equation}
Applying \eqref{eq_admm_iter} on \eqref{eq_var_opt_for_admm} and using a denoiser in place of the proximal operator, the PnP-ADMM iterations can be derived as follows, where $\rho>0$ is an appropriately chosen step size parameter: 
\begin{align*}\label{eqs:PnP_admm} \tag{PnP-ADMM}
\begin{cases}
    x_{k+1}  = \left(K^\top K + \rho\,\id\right)^{-1}\left(K^\top y+\rho\,(z_k-u_k)\right),\\
    z_{k+1} = D\left(x_{k+1}+u_k\right), \\
    u_{k+1} = u_k+\left(x_{k+1}-z_{k+1}\right).
\end{cases}
\end{align*}

\subsection{Accelerated Convex Methods}
From the previous reformulation of plug-and-play methods as non-convex optimization of some explicit functionals, a natural question is whether or not \rev{the optimization and therefore reconstruction} can be accelerated. There are two main ways of acceleration in the optimization literature, namely via momentum and preconditioning. We note that there are no available convergence guarantees in the former case of momentum-based accelerated PnP, other than spectral analyses for linear denoisers and linear inverse problems \cite{sinha2025fista}. In gradient-based optimization, preconditioning refers to multiplying the gradient by some (usually positive definite) preconditioner matrix, to make the problem ``less ill-conditioned'' and achieve faster convergence using larger step sizes. Common instances include Newton's method, Riemannian gradient descent, or the more exotic mirror descent. In the proximal splitting case, however, the preconditioning affects not only the gradient step but also the proximal step. For example, the preconditioned proximal gradient method to minimize $f+g$ takes the form
\begin{subequations}
\begin{gather}
    x_{k+1} =  \prox_g^{B_k} (\rev{x_k}- B_k^{-1} \nabla f(x_k)),\\
    \prox_g^{B_k}(x) = \argmin_{y \in \R^n} g(y) + \frac{1}{2} (y-x)^\top B_k (y-x).
\end{gather}
\end{subequations}
For Newton-type acceleration, the matrix $B_k$ would be replaced by the Hessians $B_k = \nabla^2 f(x_k)$, or some approximation thereof for quasi-Newton methods. The variable preconditioning in the proximal precludes the direct use of a single denoiser to replace the proximal step. 

In \cite{hong2024provable}, the authors consider both \rev{fixed and iteration-dependent} preconditioners, with applications to MRI reconstruction. For the \rev{iteration-dependent} preconditioning, convergence of the variable proximal scaling assumes a ``normalization-equivariant denoiser'' $D$, which is assumed to have the property that for any $\mu >0$ and $\Delta \in \mathbb{C}$, that $D (\mu x + \Delta \mathbf{1}) = \mu D(x) + \Delta \mathbf{1}$ \rev{where $\mathbf{1}$ represents the vector of all ones}. \rev{Denoisers can also be made to be ``adjustable'', by allowing them to take an additional input corresponding to the preconditioner. An application with diagonal preconditioners} to PnP-ADMM is considered in \cite{pendu2023preconditioned}, with applications to Poisson denoising. 

The work in \cite{tan2024provably} introduces PnP-LBFGS as a method of Newton-type acceleration that entirely bypasses the inclusion of the preconditioner in the proximal step, based on the Minimizing Forward-Backward Envelope (MINFBE) algorithm \cite{stella2017forward}. They theoretically show superlinear convergence to fixed points of a non-convex functional under standard quasi-Newton assumptions, which leads to significant empirical accelerations.


\section{Learning Image Priors using Denoisers}
In this section, we will discuss how denoisers can be used to construct an explicit prior, in contrast with proximal PnP schemes where denoisers provide implicit regularization. \rev{In the case where denoisers are used to target Gaussian noise, henceforth known as Gaussian image denoisers, a theoretical link towards Bayesian inference can be drawn using Tweedie's formula \cite{tweedie_main}.}  
\subsection{Tweedie's Formula}
\label{sec_tweedie}
Consider the task of estimating the clean image $x$ with a probability density function (p.d.f.) $p(x)$ from its noisy measurement $x_{\sigma}=x+\sigma\, w$, where $\sigma>0$ is the noise standard deviation and $w\sim\mathcal{N}(0,\text{Id})$. The p.d.f. of the noisy image $x_{\sigma}$ is given by
\begin{equation}
   p_{\sigma}(x_{\sigma}) = \int_{\mathbb{R}^n}  p(x_{\sigma}|x)\,p(x)\,\mathrm{d}x,
   \label{eq_noisy_marginal}
\end{equation}
where $\displaystyle p(x_{\sigma}|x) = \frac{1}{\left(\sigma\sqrt{2\pi}\right)^n}\exp\left(-\frac{1}{2\sigma^2}\|x_{\sigma}-x\|_2^2\right)$ is the conditional density of $x_{\sigma}$ given the clean image $x$. Differentiating $p(x_{\sigma}|x)$ with respect to $x_{\sigma}$ yields
\begin{align}
    \nabla_{x_{\sigma}}\, p(x_{\sigma}|x) &= \frac{1}{\left(\sigma\sqrt{2\pi}\right)^n}\exp\left(-\frac{1}{2\sigma^2}\|x_{\sigma}-x\|_2^2\right)\cdot \left(-\frac{1}{\sigma^2}(x_{\sigma}-x)\right)\nonumber\\
    & = -\frac{1}{\sigma^2}(x_{\sigma}-x) \cdot p(x_{\sigma}|x).
    \label{eq_noisy_conditional}
\end{align}
Now, differentiating both sides of \eqref{eq_noisy_marginal} and using the identity in \eqref{eq_noisy_conditional}, we get
\begin{align}
   \nabla_{x_{\sigma}}\,p_{\sigma}(x_{\sigma}) &= \int_{\mathbb{R}^n}  \left(\nabla_{x_{\sigma}}\,p(x_{\sigma}|x)\right)\,p(x)\,\mathrm{d}x= -\frac{1}{\sigma^2}\int_{\mathbb{R}^n}(x_{\sigma}-x)p(x_{\sigma}|x)\,\,p(x)\,\mathrm{d}x\nonumber\\
   & = -\frac{x_{\sigma}}{\sigma^2}\int_{\mathbb{R}^n} p(x_{\sigma}|x)\,\,p(x)\,\mathrm{d}x+\frac{1}{\sigma^2}\int_{\mathbb{R}^n}x\,p(x_{\sigma}|x)\,\,p(x)\,\mathrm{d}x\nonumber\\
   &= -\frac{x_{\sigma}}{\sigma^2}\,\,p_{\sigma}(x_{\sigma})+\frac{1}{\sigma^2}\int_{\mathbb{R}^n}x\,p(x_{\sigma}|x)\,\,p(x)\,\mathrm{d}x.
   \label{eq_grad_noisy}
\end{align}
Subsequently, dividing both sides of \eqref{eq_grad_noisy} by $p_{\sigma}(x_{\sigma})$ leads to
\begin{align}
    \frac{\nabla_{x_{\sigma}}\,p_{\sigma}(x_{\sigma})}{p_{\sigma}(x_{\sigma})} =-\frac{x_{\sigma}}{\sigma^2} + \frac{1}{\sigma^2}\int_{\mathbb{R}^n}x\,\, \left(\frac{p(x_{\sigma}|x)\,\,p(x)}{p_{\sigma}(x_{\sigma})}\right)\,\,\mathrm{d}x.
    \label{eq_grad_noisy_1}
\end{align}
Since $\displaystyle \frac{\nabla_{x_{\sigma}}\,p_{\sigma}(x_{\sigma})}{p_{\sigma}(x_{\sigma})} =\nabla_{x_{\sigma}} \log p_{\sigma}(x_{\sigma})$ and $\displaystyle \frac{p(x_{\sigma}|x)\,\,p(x)}{p_{\sigma}(x_{\sigma})}=p(x|x_{\sigma})$, thanks to Bayes rule, one can write \eqref{eq_grad_noisy_1} as
\begin{equation}
    \nabla_{x_{\sigma}} \log p_{\sigma}(x_{\sigma}) = -\frac{x_{\sigma}}{\sigma^2}+\frac{1}{\sigma^2}\int_{\mathbb{R}^n}x\,\,p(x|x_{\sigma}) \,\,\mathrm{d}x=-\frac{x_{\sigma}}{\sigma^2}+\frac{1}{\sigma^2}\,\mathbb{E}\left[x|x_{\sigma}\right].
    \label{eq_grad_noisy_2}
\end{equation}
Rearranging the terms in \eqref{eq_grad_noisy_2} leads to the familiar Tweedie's formula, a direct connection between the optimal \rev{minimum mean squared error (MMSE)} Gaussian denoiser \rev{$\mathbb{E}[x|x_\sigma]$}, and the \textit{score function} $\nabla_{x_{\sigma}} \log p_{\sigma}(x_{\sigma})$ of the noisy image $x_{\sigma}$:
\begin{equation}
        \mathbb{E}\left[x|x_{\sigma}\right] -x_{\sigma} = \sigma^2 \,\nabla_{x_{\sigma}} \log p_{\sigma}(x_{\sigma}).
        \label{eq_tweedies_final}
\end{equation}
The optimal Gaussian image denoiser essentially seeks to approximate the conditional mean of the clean image given its noisy observation (as given in \eqref{eq_tweedies_final}), and hence provides a way for approximating the noisy score function. Such a denoiser is learned by minimizing (an empirical approximation of) the MSE loss $\mathbb{E}_{x,x_{\sigma}}\left\|D(x_{\sigma})-x\right\|_2^2$ on a training dataset having clean and noisy image pairs.    

\subsection{Regularization-by-Denoising (RED)}
\label{sec_RED}
Tweedie's formula is inherently related to the Regularization-by-Denoising (RED) scheme \cite{romano2017little,reehorst2019regularization}, a variant of plug-and-play that leverages a denoiser to approximate the score function. Suppose we are given a denoiser (typically a data-driven one) which approximates the posterior mean $\mathbb{E}\left[x|x_{\sigma}\right]$. One can then construct the RED term by noting the following approximation:
\begin{equation}\label{eq:TweedieApproxRED}
    x_{\sigma} - D(x_{\sigma}) \approx -\sigma^2 \,\nabla_{x_{\sigma}} \log p_{\sigma}(x_{\sigma}),
\end{equation}
i.e. the term $x_{\sigma} - D(x_{\sigma})$ approximates the negative of the score function. This term can be immediately plugged into a gradient step similar to PnP:
\begin{equation}\label{eq:REDAdj}
    x_{k+1} = x_k - \eta [K^\top (Kx_k - y) + \frac{\lambda}{\sigma^2}(x_k - D(x_k))].
\end{equation}
RED was first proposed by Romano et al \cite{romano2017little}. The original motivation was to build an explicit regularization $R(x)$ using a denoiser function $D(x)$:
\begin{equation}
    R(x) \coloneqq x^\top(x- D(x)),
\end{equation}
with the hope that the gradient of $R(x)$ is simply $x - D(x)$. To make this true, the denoiser needs to satisfy the local-homogeneity condition:
\begin{equation}
    (1+\delta)D(x) = D((1+\delta)x), \quad \text{\,\,for all\,\,}x,
\end{equation}
and sufficiently small $\delta \in \R\setminus \{0\}$, as well as symmetry of the denoiser's Jacobian:
\begin{equation}
    JD(x)^\top= JD(x).
\end{equation}
However, Reehorst and Schniter \cite{reehorst2019regularization} later clarified that most real-world denoisers do not satisfy the Jacobian symmetry condition; hence, this view of RED is incorrect. The true gradient of $R(x)$ is instead (see \cite[Lem. 2]{reehorst2019regularization}):
\begin{equation}
    \nabla R(x) = x - \frac{1}{2}D(x) - \frac{1}{2} JD(x)^\top x 
\end{equation}
when $D$ has a nonsymmetric Jacobian, which is the case for both non-local filters (e.g., NLM, BM3D, and TNRD \cite{tnrd_pock}) and deep denoisers (such as DnCNN). If the denoiser’s Jacobian is not symmetric, then a remarkable result shows that we cannot construct any explicit regularizer whose gradient has exactly the desired form of the denoising residual $x - D(x)$.
\begin{theorem}
    [{Impossibility of explicit regularization \cite[Thm. 1]{reehorst2019regularization}}] If the denoiser $D$ has \rev{an} asymmetric Jacobian, then there is no regularization $R(x)$ that satisfies $\nabla R(x) = x- D(x)$.
\end{theorem}
\subsection{Convergence Theory for RED}
Since it is impossible to find an explicit regularizer which exactly satisfies $\nabla R(x) = x- D(x)$, Reehorst and Schniter \cite{reehorst2019regularization} consider RED's convergence to a fixed point $x^*$, satisfying
\begin{equation}\label{fc}
 K^\top (K x^* - y) + \lambda \left(x^* - D(x^*)\right) = 0.
\end{equation}
To demonstrate \rev{convergence}, we consider a provably convergent RED algorithm, namely proximal gradient RED (RED-PG) \cite{reehorst2019regularization}. Given a data fidelity $f(x)$ such as the least-squares error $\|Kx -y\|_2^2$, parameters \rev{$\lambda > 0$, $L > 0$, and initialization $v_0$}, the iterations of RED-PG are described as follows:
\begin{equation} \label{eqs:RED-PG} \tag{RED-PG}
    \begin{cases}
        x_k = \arg\min_x \left\{ f(x) + \frac{\lambda L}{2} \|x - v_{k-1}\|^2 \right\},\\
       v_k = \frac{1}{L} D(x_k) - \frac{1 - L}{L} x_k.
    \end{cases}
\end{equation}
The basic RED-PG iteration can alternatively be written as iterating the operator \( T(x) \), defined by:
\begin{equation}
T(x) \coloneqq \argmin_{z} \left\{ f(z) + \frac{\lambda L}{2} \left\| z - \left( \frac{1}{L} D(x) - \frac{1 - L}{L} x \right) \right\|^2 \right\}.
\label{eq:T_operator}
\end{equation}
This can be equivalently written using the proximal operator as
\begin{equation}
T(x) = \prox_{f/ (\lambda L)} \left( \frac{1}{L} \big( D(x) - (1 - L) x \big) \right).
\label{eq_red_alt_prox}
\end{equation}
Using the link between Tweedie's formula and RED \labelcref{eq:TweedieApproxRED}, \rev{the argument of the proximal} may be interpreted as an approximate gradient ascent step on the log-prior, since \[v_k = \frac{1}{L} D(x_k) - \frac{1 - L}{L} x_k  = x_k - \frac{1}{L}(x_k - D(x_k)) \approx x_k + \frac{1}{L}
\nabla \log p_\sigma(x_k),\] which is followed by the proximal step on the data fit $f$. 
For the RED-PG algorithm, it is easy to prove that the operator $T$ is $\alpha$-averaged, that is, $T(x) = \alpha M(x) + (1-\alpha)T(x)$ for some non-expansive operator $M$:
\begin{lemma}
[{\cite[Lem. 5]{reehorst2019regularization}}] If \( f(\cdot) \) is proper, convex, and continuous; \( D(\cdot) \) is non-expansive; and \( L > 1 \), then the operator \( T(\cdot) \) defined in \eqref{eq:T_operator} is \( \alpha \)-averaged with $\displaystyle \alpha = \max\left\{ \frac{2}{1 + L}, \frac{2}{3} \right\}$.
\end{lemma}
\noindent Using this, convergence of RED-PG to the fixed point can be proven for non-expansive denoisers, as stated below. 
\begin{theorem}[{\cite[Thm. 2]{reehorst2019regularization}}]
If \( f(\cdot) \) is proper, convex, and continuous; \( D(\cdot) \) is non-expansive; \( L > 1 \); and the operator \( T(\cdot) \) defined in \eqref{eq:T_operator} has at least one fixed point, then the RED-PG algorithm converges.
\end{theorem}
\noindent Such an approximate gradient descent step can also be accelerated via Nesterov's momentum as described in \cite[Alg. 6]{reehorst2019regularization}: 
\begin{equation} \label{eqs:RED-APG} \tag{RED-APG}
    \begin{cases}
        x_k = \arg\min_x \left\{ f(x) + \frac{\lambda L}{2} \|x - v_{k-1}\|^2 \right\},\\
        t_k = \frac{1+\sqrt{1+4t_{k-1}^2}}{2},\\
        z_k = x_k + \frac{t_{k-1} -1}{t_k}(x_k - x_{k-1}),\\
       v_k = \frac{1}{L} D(x_k) - \frac{1 - L}{L} z_k.
    \end{cases}
\end{equation}
Similarly to the lack of convergence results for momentum-based PnP, theoretical convergence of the RED-APG is not yet established and is still an open question. 








\section{Convergence of Proximal PnP Methods}
For stability and mathematical interpretation, it is important to guarantee the convergence of PnP iterations. Since (proximal) PnP methods are derived by replacing proximal operators using off-the-shelf denoisers, their convergence is not automatically guaranteed with a generic denoiser. In this section, we will define some key notions of convergence for PnP methods from weak to strong, and highlight their practical significance. We will also present some recent representative foundational convergence theorems under each category\rev{; for additional results under each category, we refer to \Cref{tab:convergence} and towards the references}.
\subsection{Fixed-Point/Iterate Convergence}
This is by far the weakest form of convergence, which requires only that the PnP \rev{iterations converge} to a solution \rev{of some fixed point problem}. More formally, the studies on fixed-point convergence consider the PnP iterations as a fixed-point update rule of the form $x_{k+1}=\Op{T}(x_k)$, and seek to determine whether $x_k$ converges to a fixed-point $x^*$ of $\Op{T}$. The specific structure of the operator $\Op{T}$ depends on the choice of the proximal splitting algorithm, the forward operator, and the denoiser. For instance, for PnP-PGD, the iterations are given by $x_{k+1}=D\left(x_k-\eta\,\nabla f(x_k)\right)$, resulting in $\Op{T}=D\circ \left(\id-\eta\nabla f\right)$. We say that a PnP method is fixed-point convergent if $\Op{T}$ has a unique fixed-point $x^*$ (i.e., $\Op{T}(x^*)=x^*$) and the PnP iterations converge to $x^*$, meaning that $\underset{k\rightarrow\infty}{\lim}\,\,x_k=x^*$. Achieving fixed-point convergence of a PnP algorithm essentially boils down to ensuring that $\Op{T}$ is a contraction mapping (under suitable conditions on the forward operator and the denoiser). This mode of convergence inherently guarantees that the solution does not worsen as the PnP iterations are repeated a large number of times.    

\begin{theorem}[Fixed-point convergence of PnP-DRSdiff {\cite[Thm. 3]{ryu2019plug}}]\label{thm:pnp_admm_fp}
    \leavevmode \\ Consider the PnP-DRSdiff algorithm, given by the iterative updates \rev{
	\begin{equation}\label{pnp_drs1}
    \begin{cases}
        y_{k}  =\prox_{\tau f}\left(x_k\right),\\ z_{k} = D(2y_{k}-x_k),  \\ 
        x_{k+1}=x_k+z_{k}-y_{k},
    \end{cases}
	\end{equation}}
	where the data-fidelity term $f$ is $\mu$-strongly convex. \rev{Letting $\Op{T}$ be the following operator
	\begin{align}
		\Op{T}=\frac{1}{2}\id + \frac{1}{2}\left(2D-\id\right)\left(2\prox_{\tau f}-\id\right),
		\label{pnp_drs_fp}
	\end{align}
	one may equivalently express \eqref{pnp_drs1} as a fixed-point iteration of the form $x_{k+1}=\Op{T}(x_k)$.} Suppose that the denoiser $D$ satisfies
	\begin{equation}
		\left\|\left(D-\id\right)(u)-\left(D-\id\right)(v)\right\|_2 \leq \epsilon  \left\|u-v\right\|_2, 
		\label{cond_denoiser_pnpDRS}
	\end{equation}
	for all $u,v\in \rev{\mathcal{X}}$ and some $\epsilon>0$, and the strong convexity parameter $\mu$ is such that $\displaystyle\frac{\epsilon}{(1+\epsilon-2\epsilon^2)\,\mu}<\tau$ is satisfied. Then the operator $\Op{T}$ is contractive and the \rev{PnP-DRSdiff} algorithm is fixed-point convergent. 
\end{theorem}
\begin{remark}
    As noted in \cite{ryu2019plug}, fixed-point convergence of PnP-DRSdiff follows from monotone operator theory if $\left(2D-\id\right)$ is non-expansive, but \eqref{cond_denoiser_pnpDRS} imposes a less restrictive condition on the denoiser. Additionally, the data fidelity term is assumed to be strongly convex, which does not hold for ill-posed inverse problems (when the forward operator has a non-trivial null space). 
\end{remark}
\noindent We note further that global fixed-point convergence in this sense implies that the reconstruction is independent of the initialization. In practice, the global non-expansiveness does not hold, and instead is softly enforced to hold locally. This allows for convergence from different initializations to different fixed points.
\rev{
\subsection{Kurdyka--Łojasiewicz Property}
In the absence of convexity of the prior or contractivity of the denoiser, another weaker form of convergence utilizes the Kurdyka--Łojasiewicz (KL) property of a function. This is a general property that may be satisfied using architectural choices on neural networks. In the following, $\partial^l$ denotes the limiting subdifferential.

\begin{definition}[Kurdyka--Łojasiewicz property \cite{attouch2013convergence,bolte2014proximal}]
    Let $\varphi:\R^n \rightarrow \R \cup \{+\infty\}$ be a proper and lower semi-continuous function. $\varphi$ satisfies the \emph{Kurdyka--Łojasiewicz (KL) property} at a point $x_*$ in $\dom \partial^l \varphi$ if there exists $\eta \in (0,+\infty]$, a neighbourhood $U$ of $x_*$ and a continuous concave function $\Psi:[0,\eta) \rightarrow [0,+\infty)$ such that:
    \begin{enumerate}
        \item $\Psi(0) = 0$;
        \item $\Psi$ is $\mathcal{C}^1$ on $(0,\eta)$;
        \item $\Psi'(s)>0$ for $s \in (0,\eta)$;
        \item For all $u \in U \cap \{\varphi(x_*) < \varphi(u) < \varphi(x_*) + \eta\}$, we have
        \[\Psi'(\varphi(u) - \varphi(x_*)) \dist (0, \partial^l \varphi(u)) \ge 1.\]
    \end{enumerate}
    We say that $\varphi$ is a KL function if the KL property is satisfied at every point of $\dom \partial^l \varphi$.
\end{definition}
\noindent The KL property can be interpreted as a certain regularity condition on a function $\varphi$. Indeed, if $\varphi$ satisfies the KL property at a critical point $\overline{u}$, then it can be shown that subgradients of $u\mapsto \Psi (\varphi(u) - \varphi(\overline{u}))$ have norm bounded away from one, also known as \textit{sharpness} \cite{bolte2014proximal}. This geometric property ensures that critical points have sufficient regularity properties. While the definition is seemingly complicated, the KL property holds for many classes of functions, with some examples given in the following proposition.
\begin{proposition}
    The following classes of functions satisfy the KL property \cite{attouch2009convergence,attouch2010proximal}:
    \begin{enumerate}
        \item Subanalytic functions that are continuous on their domain (including analytic functions);
        \item Uniformly convex functions $f$, for which there exists some $K>0, p\ge 1,$ such that for all $x,y\in \mathcal{X},\, u \in \partial f(x)$,
        \[f(y) \ge f(x) + \langle u, y-x \rangle + K\|y-x\|^p;\]
        \item Semialgebraic functions, which are functions whose graphs are finite unions of the form 
        \[\{x \in \R^{d+1} \mid p_i(x) = 0, q_i(x) < 0, i=1,...,p\},\]
        where $p_i, q_i$ are polynomials.
    \end{enumerate}
    In particular, the following are also semialgebraic \cite{bochnak2013real}:
    \begin{enumerate}
        \item Finite sums and products of semialgebraic functions;
        \item Compositions of semialgebraic functions or mappings;
        \item Indicator functions of semialgebraic sets; and
        \item Generalized inverses of semialgebraic mappings.
    \end{enumerate}
\end{proposition}
\noindent The (sub)analytic functions and semialgebraic characterizations are particularly useful. \rev{As a special case,} smooth neural networks and networks with piecewise-polynomial activations, such as ReLU, both satisfy the KL property. In the absence of convexity, one can \rev{instead} use the KL property to show convergence. For example, the following abstract theorem shows convergence under certain conditions on the iterates. Moreover, the convergence is fast in the sense that the sequence has finite length. 
\begin{theorem}[{\cite[Thm. 2.9]{attouch2013convergence}}]
    Let $f:\R^n \rightarrow \R \cup \{+\infty\}$ be a proper, lower semi-continuous function. Suppose a sequence $(x_k)_{k\in\mathbb{N}}$ and constants $a,b>0$ satisfy the following properties:
    \begin{enumerate}
        \item (Sufficient decrease). For $k \in \mathbb{N}$,
        \[f(x_{k+1}) + a \|x_{k+1} - x_k\|^2 \le f(x_k);\]
        \item (Relative error). For $k \in \mathbb{N}$, there exists $w_{k+1} \in \partial^l f(x_{k+1})$ such that
        \[\|w_{k+1}\| \le b \|x_{k+1} - x_k\|,\]
        \item (Continuity). There exists a subsequence $(x_{k_j})_{j \in \mathbb{N}}$ and a cluster point $\tilde{x} \in \R^n$ such that
        \[x_{k_j} \rightarrow \tilde{x} \text{ and } f(x_{k_j}) \rightarrow f(\tilde{x})\quad \text{as } j\rightarrow \infty.\]
    \end{enumerate}
    If, furthermore, $f$ has the KL property at the cluster point $\tilde{x}$, then the entire sequence $(x_k)_{k\in \mathbb{N}}$ converges to $\tilde{x}$. Moreover, $\tilde{x}$ is a critical point of $f$, and $(x_k)_{k \in \mathbb{N}}$ has finite length. 
\end{theorem}
\noindent The final property of the above theorem is the most pertinent in the context of PnP algorithms, as it demonstrates convergence to a critical point of the underlying (non-convex) functional \cite{hurault2022proximal,tan2024provably}. 
}
\subsection{Objective Convergence}
While convergence to a fixed point guarantees stability under repeated iterations, the fixed point generally does not lend itself to a variational interpretation. That is, the fixed point is not generally a minimizer or a stationary point of a variational energy function induced by the denoiser. Objective convergence of PnP draws a direct parallel between PnP and classical variational schemes by ensuring that the PnP solution is indeed a stationary point of some variational objective (which can potentially be non-convex depending on the denoiser). Such convergence can be shown by imposing special structures on the denoiser (while leveraging convergence analysis of proximal gradient descent or some variant of it in the non-convex setting).

\rev{One popular structure is to define the denoiser as a gradient-step (GS), result in a so-called GS denoiser. These denoisers take the form $D=\id-\nabla g$, where the ``potential function'' $g$ is proper, lower semi-continuous, and differentiable with an $L$-Lipschitz gradient. Using this structure, the denoiser can be substituted into the gradient step of the FBS: for a step-size $\tau>0$ and relaxation parameter $\lambda>0$,
\begin{align}
		x_{k+1} & = \prox_{\tau f}\left(x_k-\tau \,\lambda \,\nabla g(x_k)\right)\nonumber \\ &=\prox_{\tau  f}\circ \left(\tau\lambda\, D+\rev{(1-\tau\lambda) \id}\right)(x_k),
		\label{eq:pnp_gs_hqs}
	\end{align}
where $\circ$ denotes function composition. GS denoisers enjoy the following convergence.}
\begin{theorem}[Objective convergence of PnP iterations with GS denoisers \cite{gs_denoiser_hurault_2021}]
	\label{thm:pnp_admm_gsd}
	Suppose the denoiser is constructed as a \rev{GS} denoiser $D=\id-\nabla g$, where $g$ is proper, lower semi-continuous, and differentiable with an $L$-Lipschitz gradient. \rev{Suppose further that the data-fidelity} $f \colon \rev{\mathcal{X}}\to\R\cup \{+\infty\}$ is convex and lower semi-continuous. Then, the following guarantees hold for $\tau<\frac{1}{\lambda\, L}$:
	\begin{enumerate}
		\item The sequence $F(x_k)$, where $F=f+\lambda\, g$, is non-increasing and convergent.
		\item $\left\|x_{k+1}-x_k\right\|_2 \to 0$, which indicates that iterations are stable, in the sense that they do not diverge if one iterates indefinitely.
		\item All limit points of $\{x_k\}$ are stationary points of $F(x)$.
	\end{enumerate}
	Notably, the PnP iteration defined by \eqref{eq:pnp_gs_hqs} is exactly equivalent to proximal gradient descent on $f+\lambda\,g$, with a potentially non-convex $g$.
\end{theorem}
\noindent The construction of gradient-step denoisers is motivated by Tweedie's formula, which states that the optimal minimum mean squared-error (MMSE) Gaussian denoiser indeed has the form \rev{of a gradient step under some prior term} $D(x) = \id - \nabla g(x)$. Similarly to \cite{romano2017little}, \rev{the potential} $g(x)$ \rev{typically takes the form} $g(x) = \frac{1}{2}\left\|x-N(x)\right\|_2^2$, where $N(x)$ is taken to be \rev{some differentiable neural network such as DRUNet} \cite{zhang2021plug}. As noted in \cite{gs_denoiser_hurault_2021}, this specific design leads to a powerful denoiser while facilitating convergence analysis. The GS-PnP algorithm and following analysis can be seen to be similar to \rev{\ref{eq:PnP-PGD}}, with the gradient step and proximal step flipped. Some other recent PnP objective convergence results under specific technical assumptions on the denoiser can be found in \cite{kunal_9380942,tan2024provably}.  
\subsection{Convergent Regularization Using PnP}
While objective convergence ensures a one-to-one connection between PnP iterates and the minimization of a variational objective, it does not provide any guarantees about the regularizing properties of the solution that the iterates converge to. In the same spirit as classical regularization theory, it is therefore desirable to be able to control the implicit regularization \rev{induced} by the denoiser in PnP algorithms, and analyze the limiting behavior of the PnP reconstruction as the noise level and the regularization strength tend to zero. More precisely, assuming that the PnP iterations converge to a solution $\hat{x}\left(y^\delta,\lambda\right)$, where $\delta$ denotes the noise level and $\lambda$ is an explicit regularization penalty parameter associated with the denoiser, one would like to obtain appropriate selection rules for $\sigma$ and/or $\lambda$ such that $\hat{x}\left(y^\delta,\lambda\right)$ exhibits convergence akin to \eqref{eqn:classicRegCondition} in the limit as $\delta\rightarrow 0$. To the best of our knowledge, one of the first analyses of this kind was reported in \cite{ebner2022plugandplay}, and the precise convergence result is stated in Theorem \ref{thm:pnp_conv_reg_haltmeyer}.
\begin{theorem}[Convergent PnP regularization {\cite[Thm. 3.14]{ebner2022plugandplay}}]
	\label{thm:pnp_conv_reg_haltmeyer}
	Consider the PnP-PGD iterates corresponding to a quadratic fidelity term, which takes the form
	\begin{equation}
		x_{\lambda,k+1}^{\delta} = D_{\lambda}\left(x_{\lambda,k}^{\delta}-\eta\,K^\top\left(Kx_{\lambda,k}^{\delta}-y^{\delta}\right)\right),
		\label{eq:pnp_conv_reg_thm}
	\end{equation}
	where $D_{\lambda}$ is a denoiser with a tuneable regularization parameter $\lambda$. \rev{Suppose that the family of denoisers $\{D_\lambda\}_{\lambda>0}$ satisfies appropriate assumptions (see Definition 3.1 in \cite{ebner2022plugandplay} for details), in particular that they are contractive so that the PnP iterations converge.} Let $\PnP\left(\lambda,y^{\delta}\right)$ be the fixed point of the PnP iteration \eqref{eq:pnp_conv_reg_thm}. For any \rev{$y=y^0\in\mathscr{R}(K)$} and any sequence $\delta_k>0$ of noise levels converging to $0$, there exists a sequence $\lambda_k$ of regularization parameters converging to $0$ such that for all $y_k$ with $\|y_k-y^0\|_2\leq \delta_k$:
	\begin{enumerate}
		\item $\PnP\left(\lambda,y^{\delta}\right)$ is continuous in $y^{\delta}$ for any $\lambda>0$.
		\item The sequence $\left(\PnP\left(\lambda_k,y_k\right)\right)_{k\in\mathbb{N}}$ has a weakly convergent subsequence.
		\item The limit of every weakly convergent subsequence of $\left(\PnP\left(\lambda_k,y_k\right)\right)_{k\in\mathbb{N}}$ is a solution of the noiseless operator equation $y^0=\rev{K}x$.
	\end{enumerate}
\end{theorem}
\noindent Establishing convergence in the sense of regularization ensures that the implicit regularization effect of the denoiser vanishes to zero as the noise level in the measurement diminishes, thereby guaranteeing that there is no over- or under-regularization. Recently, the convergent regularization property of PnP algorithms with a linear denoiser was shown in \cite{Hauptmann2024_linear_pnp}, where the regularization strength of the denoiser is controlled through a spectral filtering-based approach.  
\section{Practical Constraints and Training}
To comply with the theoretical analysis, the denoisers used in PnP-like schemes need to satisfy certain constraints. In this section, we mention some of these practical constraints and how they are enforced during training. We also demonstrate the empirical performance of some recent PnP algorithms in terms of image quality and convergence of the iterates. 
\subsection{Weakly Enforced Spectral Constraints}\label{ssec:spectral}
For convergence analysis, one key requirement is that the denoiser should take the form of a proximal step. For a gradient step denoiser $D_\sigma = \id - \nabla g_\sigma$ to be a proximal operator $D_\sigma = \prox_{\phi_\sigma}$ of some weakly convex function $\phi_\sigma$, a sufficient condition is that $\nabla g_\sigma$ is $L_\sigma$-Lipschitz for some $L_\sigma<1$ \cite{gribonval2020characterization,hurault2022proximal}. Here, we use the subscript $\sigma$ to denote specifically that the denoiser is trained to remove Gaussian noise of standard deviation $\sigma$. As \rev{enforcing the Lipschitz condition through architectural choices or otherwise is difficult}, a standard approach in practice is to \rev{penalize the network} in the loss function \rev{if the Lipschitz constant is too large}. Noting the equivalence of the Lipschitz constant and the spectral norm of $\nabla^2 g_\sigma = J(\id - D_\sigma)$, this consists of adding a spectral regularization term of the form 
\begin{equation*}
    \mathbb{E}_{x\sim p, \xi \sim \mathcal{N}(0, \sigma^2)}\max(\|J(\id - D_\sigma)(x+\xi)\|, 1-\varepsilon).
\end{equation*}
\rev{Here, $\varepsilon \in (0,1)$ is a tuneable hyperparameter to control how strongly the spectral constraint should be enforced.} This penalizes the spectral norm of $\nabla^2 g_\sigma$, typically approximated using a power iteration. \rev{This method can be extended also to learn monotone operators \cite{pesquet2021learning,belkouchi2025learning}, while other} methods of softly enforcing the Lipschitz constant include (approximate) layer-wise projections onto the Stiefel manifold of orthogonal matrices \cite{anil2019sorting,cisse2017parseval}. However, as the Lipschitz constant is not strictly enforced to be less than one, the algorithms suffer from occasional divergence. 

\subsection{Backtracking for Lipschitz Control}
As mentioned in \Cref{thm:pnp_admm_gsd}, the gradient-step paradigm for PnP instead replaces the gradient step in a splitting with a denoiser, and applies the proximal operator on the fidelity term \cite{gs_denoiser_hurault_2021}. \rev{In this case, the theoretically convergent sequence $F(x_k)$ requires the computation of $F = f+\lambda g$. This is computable since $f$ is a known fidelity term, and $g$ takes the special form $g(x) = \frac{1}{2}\|x-N(x)\|$ for a neural network $N(x)$.} In this case, since the step size in the splitting is allowed to be variable, it remains to find an upper bound on the Lipschitz constant of $\nabla g_\sigma$, such that $D_{\lambda,\sigma} = \id - \lambda \nabla g_\sigma$ is a descent step. 

Instead of approximating the Lipschitz constant to find a (possibly small) appropriate step size, one may instead directly consider the consequential descent condition. This problem takes the following form: find a $\lambda \in (0,1/2)$ such that
\begin{equation*}
    F(x_k) - F(\text{GS-PnP}_{\lambda}(x_k)) \ge  \lambda^{-1} \|\text{GS-PnP}_{\lambda}(x_k) - x_k\|^2.
\end{equation*}
This can be executed similarly to an Armijo line search, and can be shown to converge in finitely many iterations under the standard assumptions.

\begin{table}[]
\centering
\caption{Summary of the properties of some convergent methods. By \textit{iterate convergence}, we mean that the entire sequence converges to a point. For methods with residual convergence, they consider convergence of the form $\min_{l \le k} \|x_{l+1}-x_l\|$. The KL property is used to transform the convergence of residuals to the convergence of iterates. \rev{Objective convergence denotes whether or not the cluster points are critical points of some computable function.} \rev{The denoisers in the latter six methods use denoisers in gradient-step form} $D_\sigma = \id - \nabla g_\sigma$, \rev{and} we denote by $L_g$ the Lipschitz constant of $g_\sigma$.}
\label{tab:convergence}
\resizebox{\textwidth}{!}{%
\begin{tabular}{@{}ccccccccc@{}}
\toprule
\multirow{2}{*}{PnP method} &
  \multirow{2}{*}{Splitting} &
  \multirow{2}{*}{Convergent?} &
  \multicolumn{3}{c}{\textbf{Denoiser constraint}} &
  \multicolumn{3}{c}{\textbf{Notion of convergence}} \\
 &
   &
   &
  Spectral &
  Line-search &
  KL &
  Iterate &
  Residual &
  Objective \\ \midrule
DPIR \cite{zhang2021plug}&
  HQS &
  \xmark &
  \multicolumn{6}{c}{\xmark} \\
PnP-ADMM \cite{chan2017plug} &
  ADMM &
  \cmark &
  \multicolumn{3}{c}{"Bounded denoiser" $\|D_\sigma x - x\|^2 \le C \sigma^2$} &
  \cmark &
  \xmark &
  \xmark \\
GS-PnP \cite{gs_denoiser_hurault_2021} &
  PGD &
  \cmark &
  $L_g<1$ &
  \cmark &
  \xmark &
  \cmark &
  $\mathcal{O}(1/\sqrt{k})$ &
  \cmark \\
PnP-PGD &
  PGD &
  \cmark &
  $L_g<1$ &
  \xmark &
   \cmark &
   \cmark &
  $\mathcal{O}(1/\sqrt{k})$ &
  \cmark \\
PnP-$\alpha$PGD &
  PGD &
  \cmark &
  $L_g<1$ &
  \xmark &
  \xmark &
  \xmark &
  $\mathcal{O}(1/\sqrt{k})$ &
  \cmark \\
PnP-DRS &
  DRS &
  \cmark  &
  $L_g<1$ &
  \xmark &
   \cmark &
   \cmark &
  $\mathcal{O}(1/\sqrt{k})$ &
  \cmark \\
PnP-DRSdiff&
  DRS &
  \cmark  &
  $L_g<1/2$ &
  \xmark &
   \cmark &
   \cmark &
  $\mathcal{O}(1/\sqrt{k})$ &
  \cmark \\
PnP-LBFGS \cite{tan2024provably}&
  MINFBE \cite{stella2017forward}&
  \cmark &
  $L_g<1$ &
  \cmark &
   \cmark &
   Superlinear &
  $\mathcal{O}(1/\sqrt{k})$ &
  \cmark \\ \bottomrule
\end{tabular}%
}
\end{table}

In Table \ref{tab:convergence}, we summarize the properties required for some recent provable PnP methods based on operator splitting convergence. \rev{To verify the convergence of the various provable PnP algorithms, we test them on a natural image deblurring task. We compare DPIR along with the latter five provable PnP algorithms in the table, where the denoiser is given by a gradient step denoiser $D_\sigma = \id - \nabla g_\sigma$, where $g_\sigma(x) = \frac{1}{2}\|x-N_\sigma(x)\|^2$ is a pretrained DRUNet architecture \cite{zhang2021plug}. To (approximately) satisfy the assumptions of the previous theorems, the Lipschitz constant of $\nabla g_\sigma$ is penalized to be less than 1 as in Section \ref{ssec:spectral}, and the activation functions are taken to be $\mathcal{C}^2$ and such that the neural network satisfies the KL property.} 

\begin{figure}[htp]%
    \centering
    \subfloat[\centering Ground Truth ]{{\includegraphics[height=2.9cm]{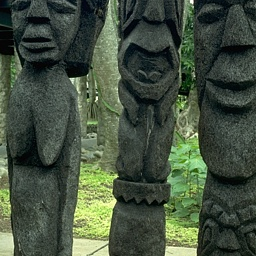}} } %
    \subfloat[\centering DPIR (24.09dB)]{{\includegraphics[height=2.9cm]{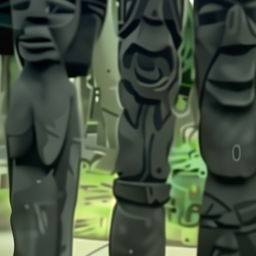}} } %
    \subfloat[\centering PnP-LBFGS (24.14dB)]{{\includegraphics[height=2.9cm]{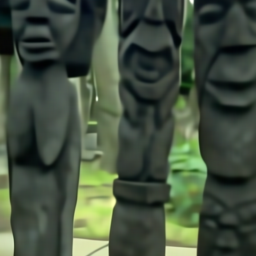}} } %
    \subfloat[\centering PnP-$\hat\alpha$PGD (24.33dB)]{{\includegraphics[height=2.9cm]{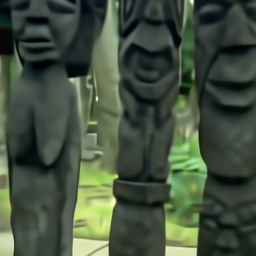}} } %

    \subfloat[\centering Corrupted ]{{\includegraphics[height=2.9cm]{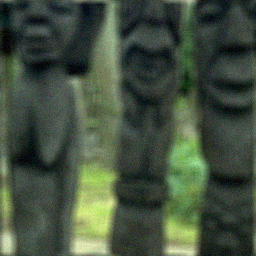}} } %
    \subfloat[\centering PnP-PGD (24.26dB)]{{\includegraphics[height=2.9cm]{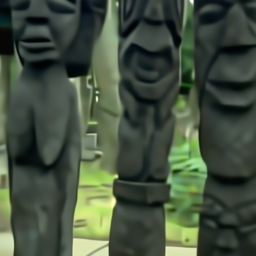}} } %
    \subfloat[\centering PnP-DRSdiff (24.26dB)]{{\includegraphics[height=2.9cm]{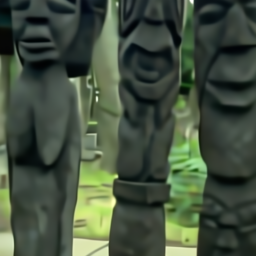}} } %
    \subfloat[\centering PnP-DRS (24.37dB)]{{\includegraphics[height=2.9cm]{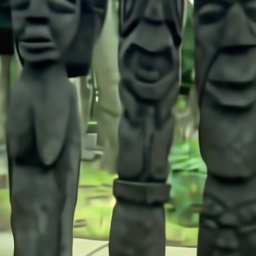}} } %
    \caption{Example reconstructions for a test image, with PSNR to ground truth in brackets. The image is blurred with a $9\times 9$ uniform blur kernel, with subsequent 3\% additive Gaussian noise. Observe that PnP-PGD and PnP-DRSdiff have the same eventual PSNR, due to targeting the same underlying functional.}
    \label{fig:deblurImages}%
\end{figure}

\begin{figure}[htp]%
    \centering
    \subfloat[\centering DPIR]{{\includegraphics[height=3cm]{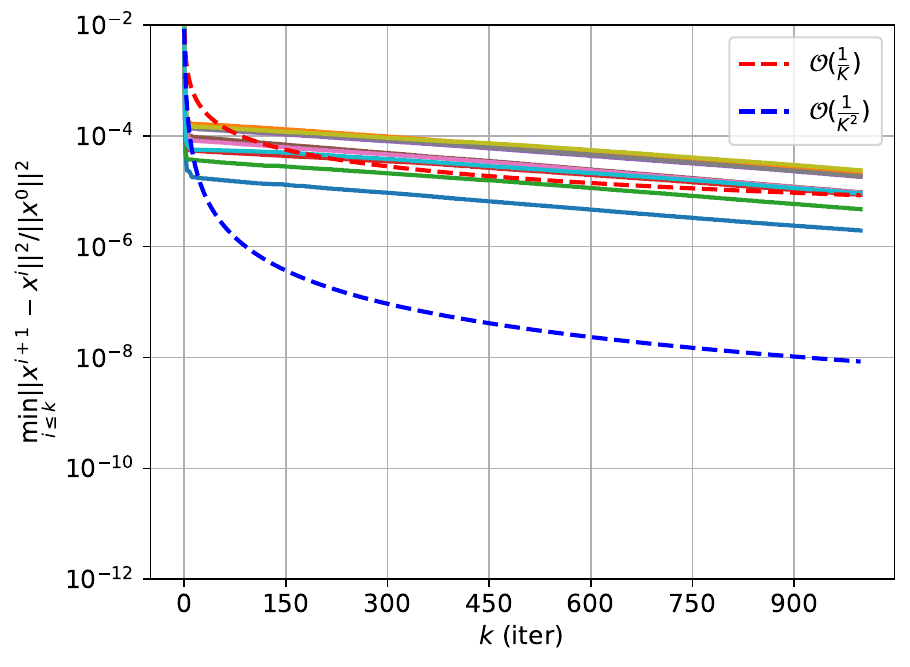}} } %
    \subfloat[\centering PnP-LBFGS]{{\includegraphics[height=3cm,trim={72px 0 0 0},clip]{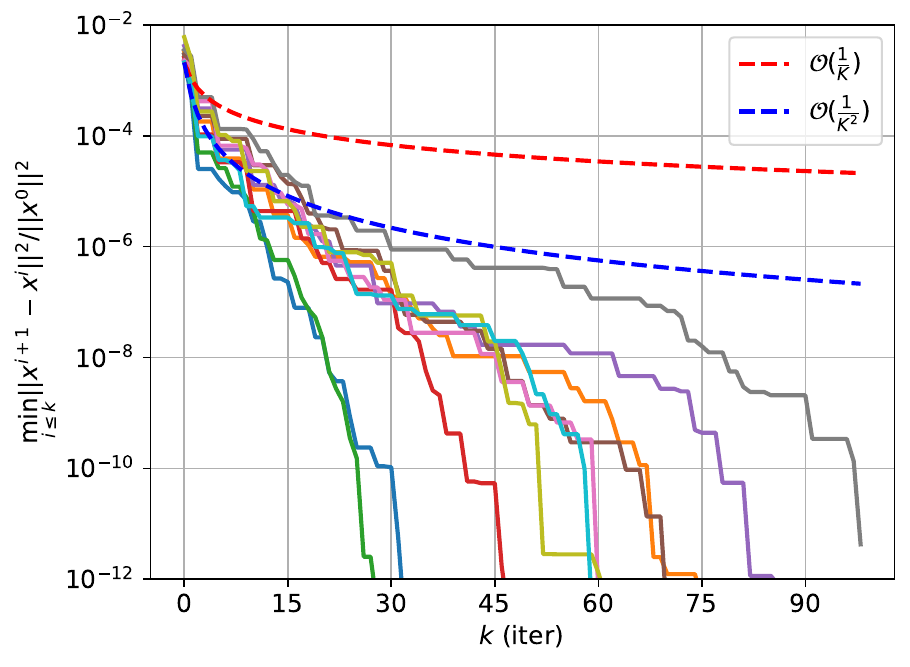}} } %
    \subfloat[\centering \rev{PnP-$\alpha$PGD}]{{\includegraphics[height=3cm,trim={72px 0 0 0},clip]{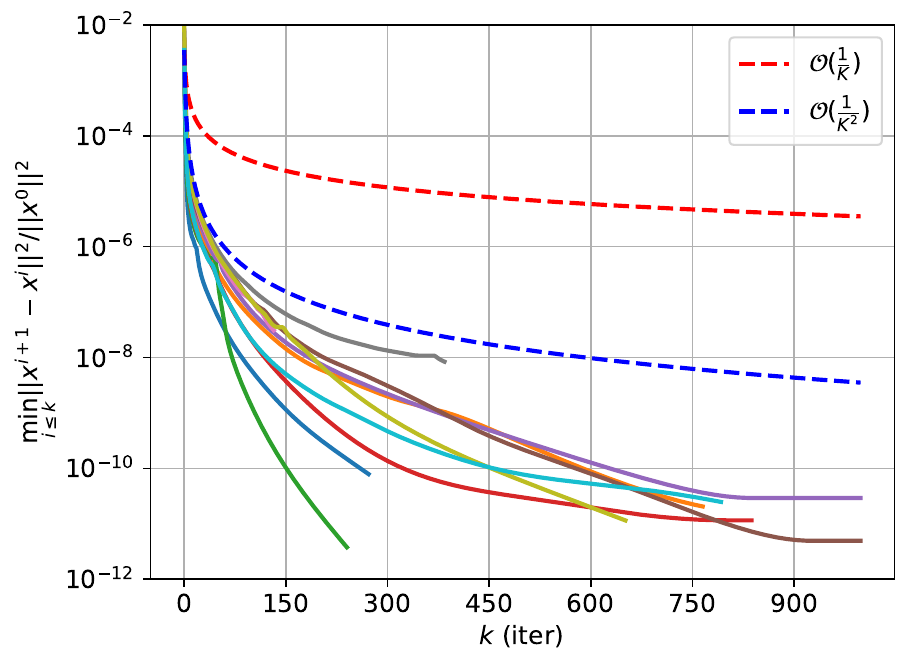}} } %

    \subfloat[\centering PnP-PGD]{{\includegraphics[height=3cm]{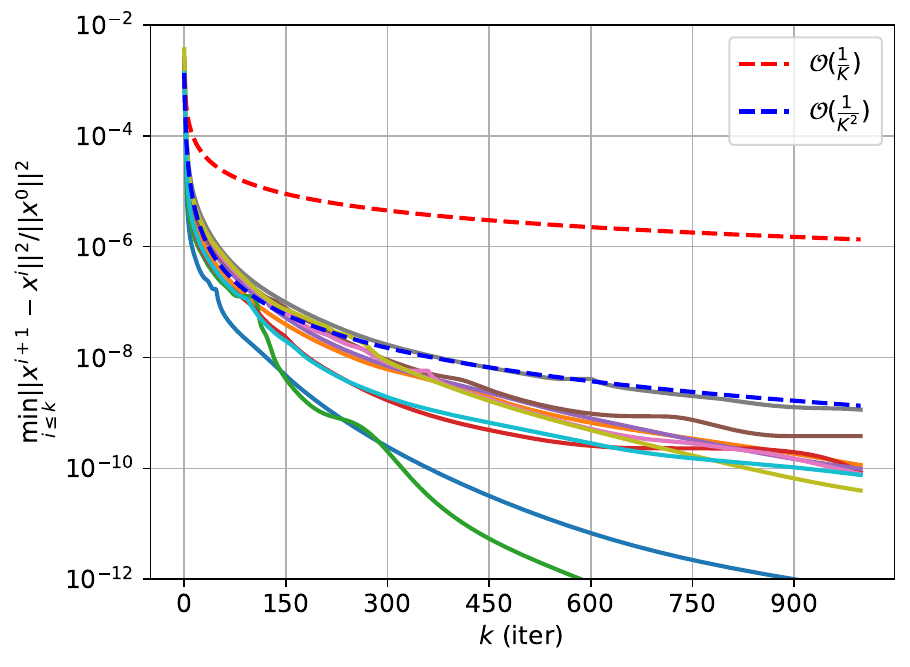}} } %
    \subfloat[\centering PnP-DRSdiff]{{\includegraphics[height=3cm,trim={72px 0 0 0},clip]{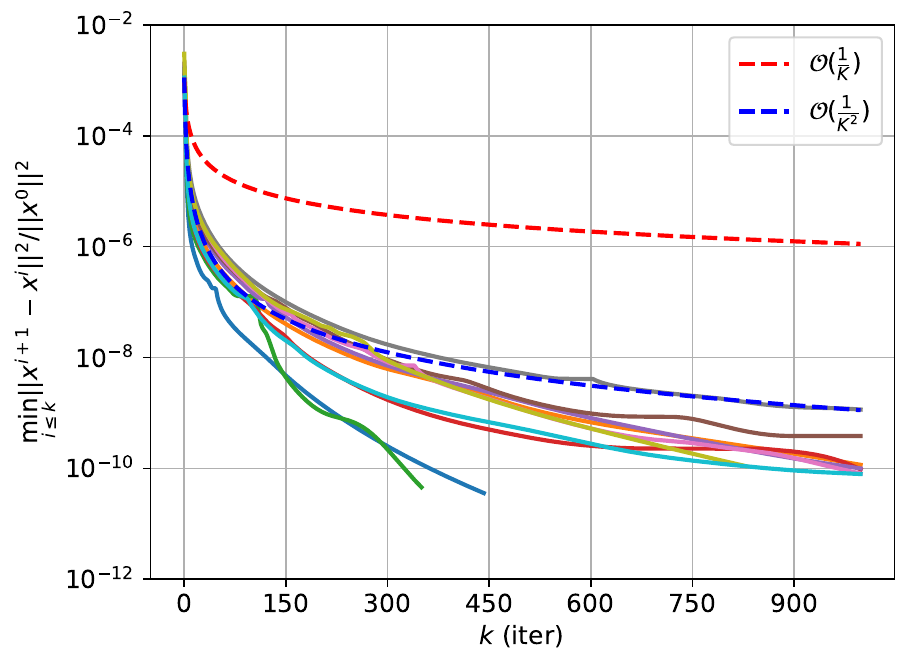}} } %
    \subfloat[\centering PnP-DRS]{{\includegraphics[height=3cm,trim={72px 0 0 0},clip]{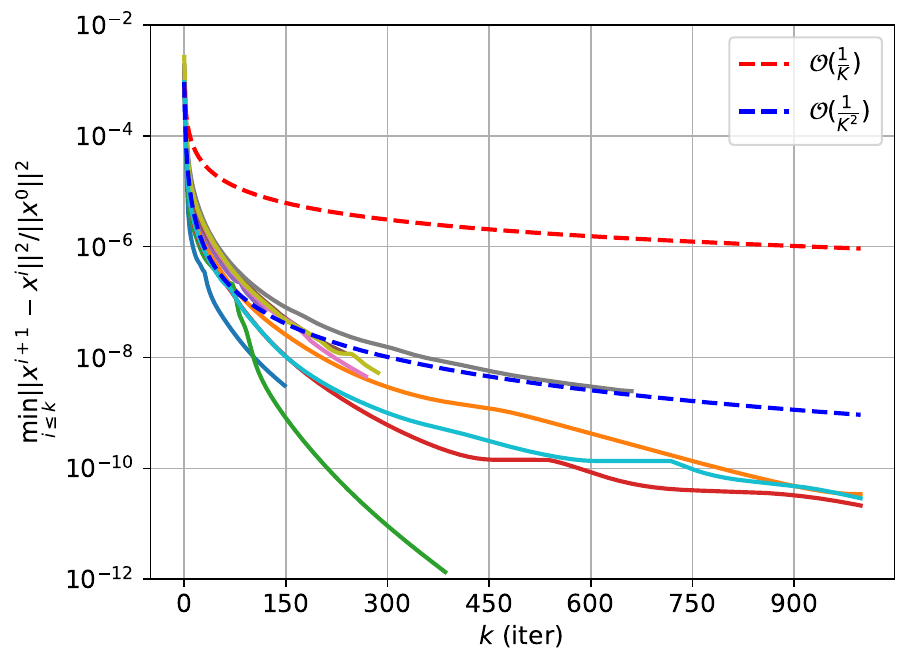}} } %
    \caption{Residual convergence for deblurring on the CBSD10 dataset, with a uniform $9\times 9$ blur kernel and 3\% additive Gaussian noise. \rev{Each solid line represents one image.} We observe that while DPIR has slow residual convergence, the provable PnP methods all have a convergent behavior, often reaching their stopping criteria given by the change in objective value. In particular, the quasi-Newton PnP-LBFGS method converges very quickly within 100 iterations.}
    \label{fig:residual}%
\end{figure}

\begin{figure}[ht]%
    \centering
    \subfloat[\centering DPIR]{{\includegraphics[height=3cm]{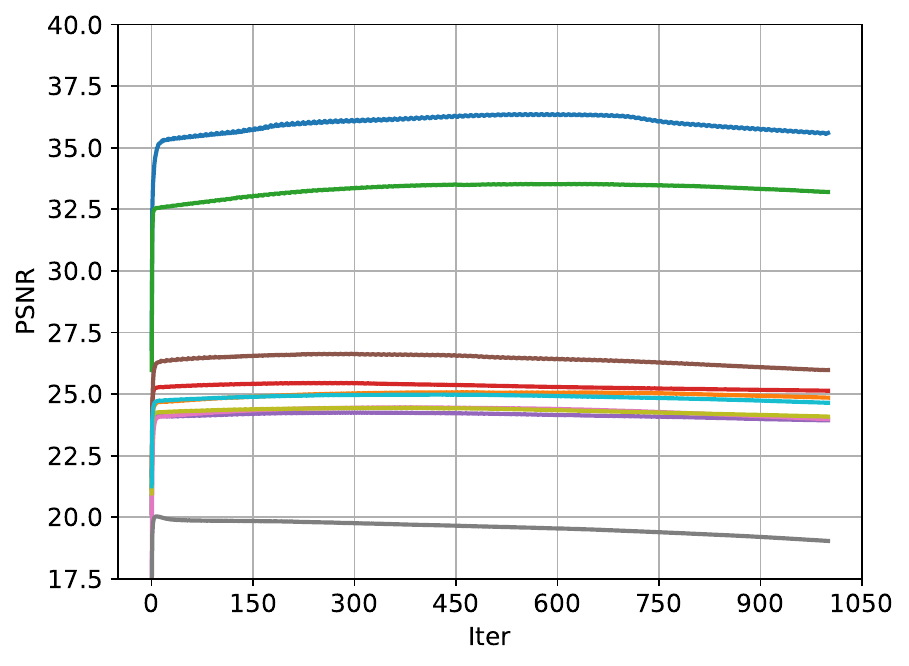}} } %
    \subfloat[\centering PnP-LBFGS]{{\includegraphics[height=3cm,trim={56px 0 0 0},clip]{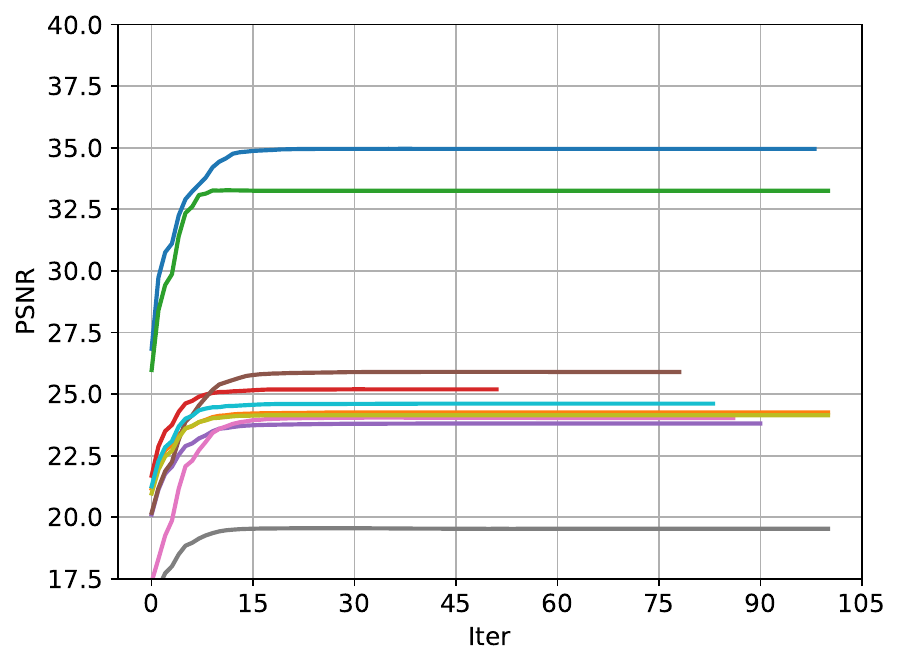}} } %
    \subfloat[\centering \rev{PnP-$\alpha$PGD}]{{\includegraphics[height=3cm,trim={56px 0 0 0},clip]{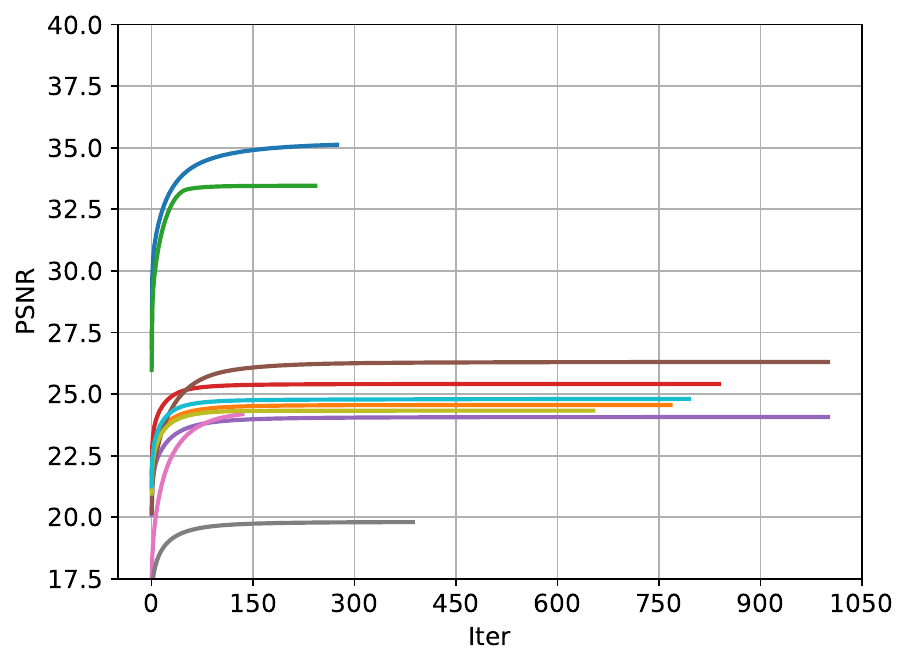}} } %

    \subfloat[\centering PnP-PGD]{{\includegraphics[height=3cm]{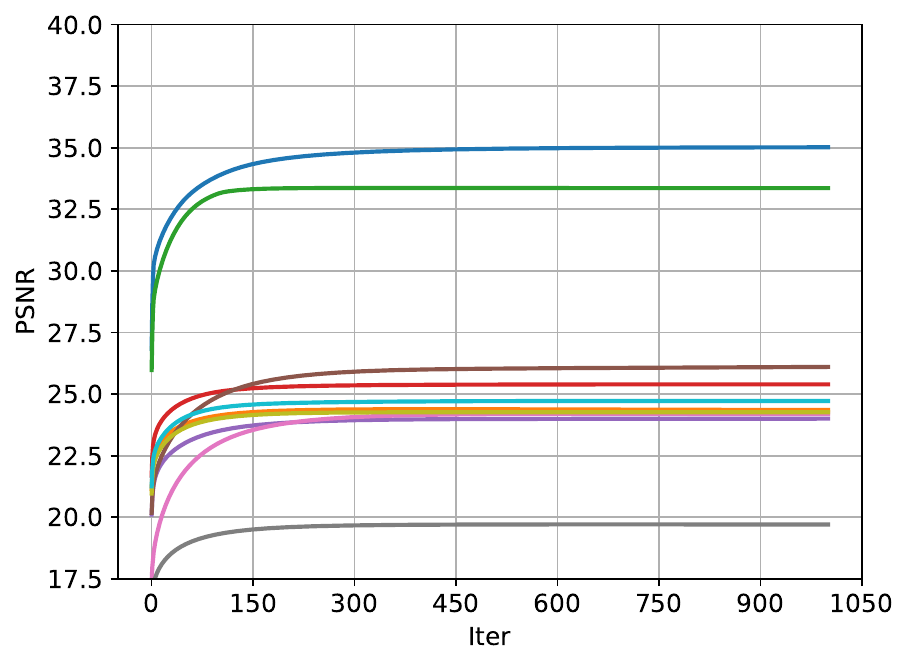}} } %
    \subfloat[\centering PnP-DRSdiff]{{\includegraphics[height=3cm,trim={56px 0 0 0},clip]{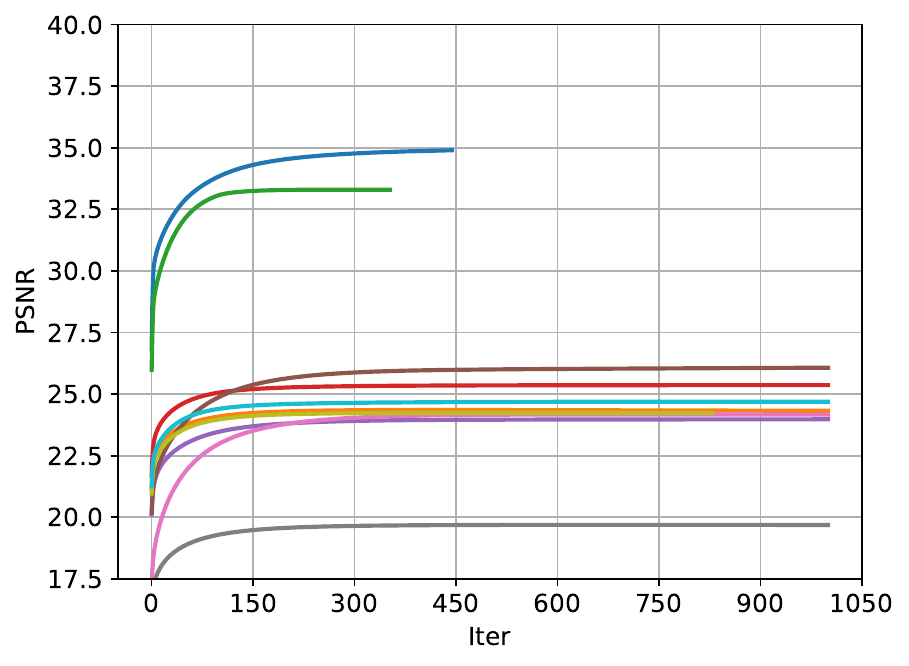}} } %
    \subfloat[\centering PnP-DRS]{{\includegraphics[height=3cm,trim={56px 0 0 0},clip]{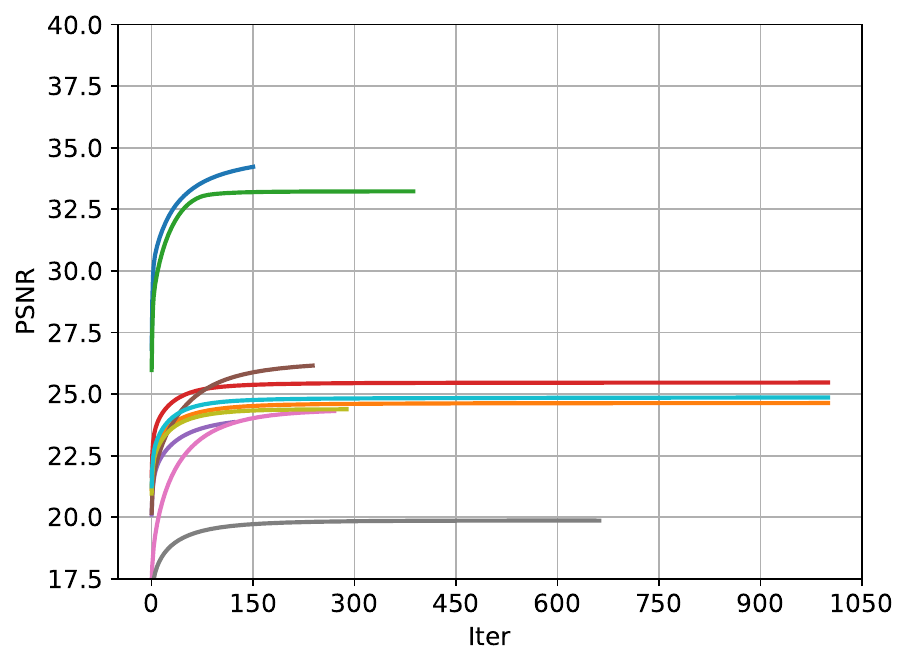}} } %
    \caption{PSNR curves for deblurring on the CBSD10 dataset, with uniform $9\times 9$ blur kernel and 3\% additive Gaussian noise. \rev{Each solid line represents one image.} We observe that the non-provable DPIR method gradually decreases in PSNR at later iterations, eventually leading to instability. In contrast, the provable PnP methods all have stable convergence curves, reaching their stopping criteria.}
    \label{fig:psnr}%
\end{figure}

Since the proof structure is quite similar for each of the methods given in Section \ref{sec1}, the constraints on the denoisers are also quite similar, and they even converge to critical points of the same functional. This is demonstrated in Figure \ref{fig:deblurImages}, where for the deconvolution task with a fixed blur kernel and PnP denoiser, the reconstructions are all fairly similar. Figures \ref{fig:residual} and \ref{fig:psnr} demonstrate the residual and peak signal-to-noise ratio (PSNR) convergence for a set of 10 images, again for the image deconvolution task. We observe that, as the theory suggests, the provable PnP methods exhibit decaying residuals and stable PSNR figures, whereas the non-provable DPIR method \cite{zhang2021plug} does not demonstrate such convergence, while deteriorating in quality as the iterations continue.

\section{Denoisers for Posterior Sampling}

The use of denoisers as components in posterior sampling has gained considerable traction as the demand for uncertainty quantification in imaging grows. In many ill-posed inverse problems, the solution space is inherently ambiguous. Instead of recovering a single best estimate, it is often more informative to approximate the full posterior distribution \( p(x | y) \). Denoisers provide a natural bridge for this, acting as powerful implicit priors that can be integrated into modern stochastic sampling schemes to draw samples from complex high-dimensional distributions.

\subsection{The Bayesian Inversion Problem}
In the Bayesian framework for inverse problems \cite{Dashti2017}, the image $x$ and the measurement $y$ are modeled as $\mathcal{X}$- and $\mathcal{Y}$-valued random variables, respectively, and the goal is to characterize the posterior distribution of $x$ given a realization of the measurement (through a summary of the posterior using a point estimate, or a mechanism that facilitates sampling from this posterior distribution, for instance). The target posterior density \( p(x| y) \) combines the likelihood \( p(y | x) \) and the image prior \( p(x) \) using Bayes' rule, modeling the image acquisition process and capturing assumptions about the clean image, respectively:
\[
p(x| y) \propto p(y |x) p(x).
\]
When the forward model is linear with additive Gaussian noise, that is, $y = Kx + w$, where $w \sim \mathcal{N}(0, \sigma_w^2 \mathbf{I})$, then the likelihood is given by $p(y | x) \propto \exp\left(-\frac{1}{2\sigma_w^2} \|y - Kx \|_2^2 \right)$, providing a link \rev{to} the fidelity term \rev{within} variational regularization. Sampling from this posterior is intractable if the prior is defined only implicitly through an image denoiser. Therefore, Monte Carlo Markov Chain (MCMC) algorithms are often used to generate samples from (an approximation to) the posterior by incorporating the image prior through an off-the-shelf pretrained image denoiser.

As discussed in Sec.~\ref{sec_RED}, modern denoisers approximate the MMSE estimate for an image corrupted by Gaussian noise. By Tweedie’s formula, the score function (gradient of the log prior) relates to the denoising operation as $\nabla_{x} \log p(x) \approx \frac{1}{\sigma^2} \big( D(x) - x\big)$, where \( D(\cdot) \) denotes the denoiser. This key observation underpins PnP-based MCMC methods and forms the basis for using denoisers within stochastic differential equation (SDE)-based generative samplers.
\subsection{Plug-and-play Denoisers for Posterior Sampling}
Score-based generative models and denoising diffusion probabilistic models (DDPMs) leverage an SDE framework to generate samples consistent with the data distribution. For the unconditional sampling, the forward SDE progressively corrupts the data with noise,
\begin{equation*}
\mathrm{d}x = f(x, t)\, \mathrm{d}t + g(t)\, \mathrm{d}w,
\end{equation*}
where \( w \) is a standard Wiener process, \rev{$f(x,t)$} defines the drift, and \rev{$g(t)$} the diffusion strength. The time-reversed SDE is given by
\[
\mathrm{d}x = \Big[ f(x, t) - g^2(t) \nabla_{x} \log p_t(x) \Big] \mathrm{d}t + g(t) \,\mathrm{d}w,
\]
which guides the generative sampling process, where the time-dependent score \( \nabla_{x} \log p_t(x) \) is approximated by a denoiser trained at varying noise scales \cite{song2021scorebased,ho2020denoising}. A particularly interesting special case is obtained by choosing \( f(x,t) \equiv 0 \) and a constant diffusion strength \rev{$g(t) \equiv 1$}. In this case, the reverse SDE simplifies to
\begin{equation}\label{eq:simplifiedRSDE}
\mathrm{d}x = -\nabla_{x} \log p_t(x)\, \mathrm{d}t + \mathrm{d}w,
\end{equation}
which is precisely the overdamped Langevin diffusion targeting the (instantaneous) data distribution with density proportional to \( p_t(x) \). Equivalently, if the denoiser provides a score surrogate via Tweedie’s formula, one may write 
\[
\mathrm{d}x \approx \left[ \rev{-\tfrac{1}{\sigma_t^2}}\left(D_{\sigma_t}(x_t) - x_t\right) \right] \mathrm{d}t + \mathrm{d}w,
\]
thereby exhibiting the plug-and-play interpretation in the unconditional setting. \rev{Here, $D_{\sigma_t}$ is a denoiser trained or calibrated for noise level $\sigma_t$, and implicitly depends on time through the noise level $\sigma_t$, which plays the role of the diffusion variance at step $t$.} A first-order Euler--Maruyama discretization with step size \( \delta>0 \) \rev{in backward time} yields the unadjusted Langevin algorithm (ULA):
\rev{
\[
x_{k+1} = x_k +\delta\, \nabla_x \log p_t(x_k) + \sqrt{\delta}\,\boldsymbol{\epsilon}_k,
\qquad \boldsymbol{\epsilon}_k \sim \mathcal{N}(0,I).
\]}
When the prior score is replaced by the denoiser-based approximation,\rev{
\[
x_{k+1} \approx x_k + \delta\, \tfrac{1}{\sigma_k^2}\big(D_{\sigma_k}(x_k) - x_k\big) + \sqrt{\delta}\,\boldsymbol{\epsilon}_k,
\]}
\noindent the iteration refines samples by alternating deterministic drift along the (approximate) score with stochastic exploration \cite{pnp_ula_marcelo}. This idea may be adapted to inverse problems \cite{kawar2022denoising}, where the goal is to sample from a distribution consistent with both the learned image prior and the measurement model (see \cite{daras2024surveydiffusionmodelsinverse,Chung2025} for recent surveys on the theory and applications of diffusion models for posterior sampling). A natural modification of the \rev{backward SDE \labelcref{eq:simplifiedRSDE}} to sample from the posterior distribution is by augmenting the drift with the likelihood score: 
\[
\mathrm{d}x =\rev{ -}\Big[ \nabla_x \log p_t(x) + \nabla_x \log p(y|x) \Big] \mathrm{d}t + \mathrm{d}w.
\]
Under the linear Gaussian model \( y = Kx + w \), \( w \sim \mathcal{N}(0,\sigma_w^2 I) \), the likelihood gradient takes the explicit form  $
\nabla_x \log p(y|x) = \tfrac{1}{\sigma_w^2} K^\top (y - Kx)$, 
so that the conditional Langevin SDE reduces to
\[
\mathrm{d}x =\rev{ -} \Big[ \nabla_x \log p_t(x) + \tfrac{1}{\sigma_w^2} K^\top (y - Kx) \Big] \mathrm{d}t + \mathrm{d}w.
\] With a denoiser-based prior score surrogate, this can be approximated as $\mathrm{d}x \approx -\Big[ \tfrac{1}{\sigma_t^2}\big(D_{\sigma_t}(x) - x\big) + \tfrac{1}{\sigma_w^2} K^\top (y - Kx) \Big] \mathrm{d}t + \mathrm{d}w$. Discretizing again via Euler--Maruyama \rev{in reverse time} gives the Langevin-type update used in practice:
\[
x_{k+1} = x_k + \delta \Big[ s_{\theta}(x_k, t_k) + \nabla_{x} \log p(y | x_k) \Big]
+ \sqrt{\delta}\, \boldsymbol{\epsilon}_k,\quad \boldsymbol{\epsilon}_k \sim \mathcal{N}(0, I).
\]
The PnP variant of the ULA scheme \cite{pnp_ula_marcelo}, known as PnP-ULA, takes the form
\[
x_{k+1} \approx x_k + \delta \Big[ \tfrac{1}{\sigma_k^2}\big(D_{\sigma_k}(x_k) - x_k\big) + \tfrac{1}{\sigma_w^2} K^\top (y - Kx_k) \Big]
+ \sqrt{\delta}\, \boldsymbol{\epsilon}_k,
\]
where the score is replaced by a Gaussian denoiser-based approximation. Here, $\sigma^2$ denotes the variance of the noise that the denoiser was trained to remove (which could be different from $\sigma_w^2$, the noise variance corrupting the measurement). In both unconditional and conditional cases, more advanced predictor--corrector strategies can be layered on top of this Euler--Maruyama backbone, with the predictor following the discretized SDE and the corrector applying a few local Langevin refinements at the same noise scale to improve stability and sampling efficiency \cite{song2021scorebased,pedrotti2024improved}. An important recent development is \emph{diffusion posterior sampling} (DPS) \cite{chung2023DPS}, which provides a principled extension of diffusion models to general noisy and nonlinear inverse problems by directly approximating posterior sampling. Unlike earlier diffusion solvers that primarily addressed noiseless linear problems, DPS incorporates the measurement model and noise statistics (e.g., Gaussian and Poisson) into the sampling dynamics through a learned time-dependent score network trained via score matching. Conceptually, DPS updates resemble PnP–ULA \cite{pnp_ula_marcelo}, in the sense that both alternate between a drift step informed by a learned prior (denoiser or score) and a stochastic step that injects noise for exploration. However, DPS departs from the strict projection-based measurement consistency by using a manifold-constrained gradient incorporated into the diffusion sampling path. This results in a stable and realistic reconstruction, particularly in challenging nonlinear and noisy inverse problems such as phase retrieval and non-uniform deblurring.


Posterior sampling with denoisers has enabled uncertainty quantification in applications such as medical image reconstruction, compressive sensing, and computational microscopy. Multiple posterior samples allow practitioners to construct pixel-wise credible intervals and detect ambiguous regions that would otherwise be hidden by deterministic estimators. However, there are several practical challenges that remain to be addressed. Convergence guarantees for these denoiser-driven SDE samplers are still limited, especially when the denoiser is highly nonlinear and trained on finite data. The computational cost of generating many samples, which often requires thousands of iterative steps, can be computationally prohibitive. \rev{Active research seeks to develop more efficient discretizations \cite{learning_disc_icml}, latent diffusion models \cite{latent_diffusion_cvpr}, or hybrid schemes that combine rapid MAP estimates with stochastic refinements \cite{NEURIPS2024_2f46ef57} to make these methods practical for large-scale problems.} Nonetheless, the use of denoisers for posterior sampling illustrates the remarkable synergy between learned priors and stochastic inference for solving high-dimensional imaging inverse problems. 
\section{Conclusions and Outlook}
In this chapter, we have surveyed the development of image denoising and the role that denoisers play in solving inverse problems through plug-and-play (PnP) methods. Beginning with classical denoising algorithms, we reviewed how modern learning-based denoisers can be seamlessly integrated into iterative schemes derived from variational regularization frameworks and proximal splitting algorithms. We discussed how PnP extends these algorithms by replacing proximal maps with powerful denoising operators, and explored related formulations such as Tweedie’s formula and the RED framework. Particular emphasis was placed on the mathematical conditions under which PnP methods converge, the constraints that must be imposed on denoisers to ensure stability, and practical considerations for deploying these techniques in real-world imaging settings.

Beyond deterministic optimization, we also briefly reviewed the use of denoisers in posterior sampling, highlighting connections to stochastic differential equations and their discretizations. This perspective bridges the gap between variational inference and generative modeling, offering a probabilistic interpretation of PnP and RED within the broader landscape of score-based methods.

Looking forward, several promising research avenues emerge. First, a deeper theoretical understanding of PnP with non-expansive yet highly expressive denoisers could relax current restrictive assumptions while retaining convergence guarantees. Second, domain-adapted and multimodal denoisers have the potential to unlock PnP applications in emerging imaging modalities and dynamic acquisition settings. Third, the intersection of PnP with diffusion-based generative priors and self-supervised learning may yield reconstruction algorithms that are simultaneously more robust and data-efficient. \rev{By embedding denoisers within SDEs and MCMC updates, one can approximate complex posteriors that would be intractable otherwise, providing both high-fidelity reconstructions and rigorous uncertainty quantification. This paradigm represents an exciting frontier for solving ill-posed inverse problems in a principled, uncertainty-aware manner.} Finally, scalable implementations capable of handling the high dimensionality and streaming nature of modern imaging data remain an important challenge, especially in time-critical domains such as medical imaging and remote sensing applications.

In summary, plug-and-play methods have evolved from an elegant algorithmic approach into a versatile and theoretically grounded framework for solving complex high-dimensional inverse problems.~With continued advances in denoiser design, theoretical analysis, and application-specific adaptation, PnP is poised to remain a central paradigm in computational imaging for years to come, with interesting theoretical and practical challenges to address.

\bibliographystyle{unsrt}
\bibliography{ref}

